\documentclass[10pt]{amsart}
\usepackage{amsmath, amsthm}
\let\oldin\in
\let\oldinfty\infty
\let\oldlangle\langle
\let\oldrangle\rangle

\usepackage{mathabx}
\renewcommand{\in}{\oldin}
\renewcommand{\infty}{\oldinfty}
\renewcommand{\langle}{\oldlangle}
\renewcommand{\rangle}{\oldrangle}

\usepackage[margin=1.5in]{geometry}
\usepackage[dvipsnames]{xcolor}
\usepackage[colorlinks, allcolors=purple!90!black]{hyperref} 
\hypersetup{
    colorlinks = true,
    linkcolor = purple!90!black,
    citecolor = Blue,
    urlcolor = Blue,
}
\usepackage{amssymb}
\usepackage[nameinlink,capitalize]{cleveref}
\usepackage{microtype} 
\usepackage{euscript}
\newcommand{\euscr}[1]{\EuScript{#1}}
\usepackage[
    backend=biber,
    style=alphabetic,
    maxbibnames=99,
    giveninits=true,
    isbn=false,
    url=false,
    doi=false,
    maxalphanames=9,
    minalphanames=4,
    backref,
]{biblatex}
\usepackage[parfill]{parskip}
\usepackage{tablefootnote}
\usepackage{thmtools}

\newcommand{\aeu}{\euscr{A}}
\newcommand{\aeudual}{\euscr{A}^\vee}

\newcommand{\upkq}{\textup{kq}}

\newcommand{\upe}{\textup{E}}
\newcommand{\uph}{\textup{H}}

\renewcommand{\ss}{\mathbb{S}}

\newcommand{\cc}{\mathbb{C}}

\newcommand{\qq}{\mathbb{Q}}

\newtheorem{thm}{Theorem}[section]

\theoremstyle{definition}

\newtheorem{example}[thm]{Example}

\theoremstyle{remark}
\newtheorem{remark}[thm]{Remark}

\numberwithin{equation}{section}

\usepackage{tikz-cd}
\usepackage{float}

\begin{document}

\title[Periodic phenomena in stable motivic homotopy theory]{Periodic phenomena in \\stable motivic homotopy theory}

\author{Jackson Morris}
\address{University of Washington, Seattle WA}
\email{\textcolor{Blue}{jacksonmorris1999@gmail.com}}


\dedicatory{To my advisors, John Palmieri and Kyle Ormsby.}


\begin{abstract}
In this survey, we study how tools from stable homotopy theory have manifested and impacted motivic homotopy theory. In particular, we discuss various motivic Adams spectral sequences, periodicity in the motivic stable homotopy groups of spheres, and synthetic spectra. We conclude with many problems for future investigation.
\end{abstract}

\maketitle

\vspace*{\fill}
\begin{figure}[H]
    \centering
    \includegraphics[scale=.3]{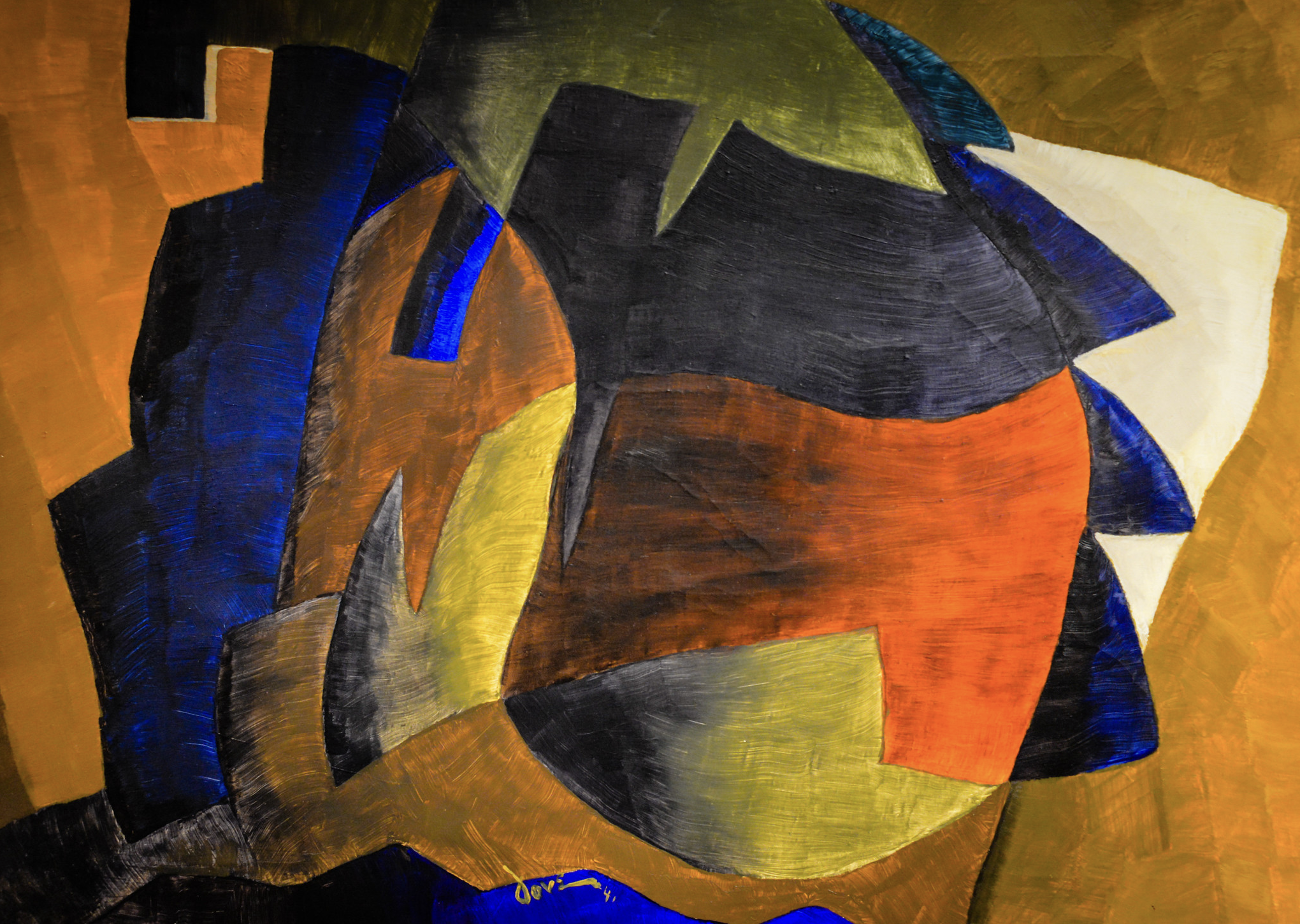}
    \caption{\emph{Morning Green}, Arthur Dove (1941).}
\end{figure}
\vspace*{\fill}
\newpage

\tableofcontents



\section{Introduction}

Motivic homotopy theory, as constructed by Morel and Voevodsky \cite{MV99}, was originally constructed to apply tools from homotopy theory to problems in algebraic geometry. Two of the motivating problems for this theory are the Milnor and Bloch--Kato conjectures, and the Lichtenbaum--Quillen conjecture \cite{MV99,Voe11}. These conjectures involve various motivic cohomology theories and their solutions are obtained by using the full power of the stable motivic homotopy category $\text{SH}(F)$. 

The objects of the stable motivic homotopy category $\mathrm{SH}(F)$ represent $\mathbb{A}^1$-invariant motivic cohomology theories on $\text{Sm}_F$ which satisfy Nisnevich descent. There is a category $\mathrm{Spc}(F)_\bullet$ of pointed motivic spaces which $\mathrm{Sm}_F$ embeds into, leading to a suspension spectrum functor
\[
\Sigma_{\mathbb{P}^1}^\infty: \mathrm{Sm}_F \to \mathrm{SH}(F).
\]
Moreover, there is a symmetric monoidal product on $\mathrm{SH}(F)$ denoted $\otimes$ whose monoidal unit is $\mathbb{S} = \Sigma_{\mathbb{P}^1}^\infty\mathrm{Spec}(F)_+$. There are a natural bigraded family of spheres 
\[
S^{s,w} = (S^1)^{\otimes s-w} \otimes (\mathbb{G}_m)^{\otimes w},
\]
leading to a bigraded family of homotopy groups $\pi_{s,w}^F(-)$. See \cite{MV99,VoeA1} for background on the foundations of motivic homotopy theory or \cite{BacElmMor25,degliseMotivicHomotopyTheory2025} for a modern approach.

The solutions to the Milnor and Bloch--Kato conjectures involves the construction of the dual motivic Steenrod algebra $\aeu^\vee = \pi_{*,*}^F(\mathrm{H}\mathbb{F}_\ell \otimes \mathrm{H}\mathbb{F}_\ell)$, the algebra of cooperations for mod-$\ell$ motivic cohomology, and interprets motivic cohomology with $\mathbb{Z}/\ell$-coefficients in terms of Milnor K-theory $\mathrm{K}^\mathrm{M}_*(F)/\ell$.  The solution to the Lichtenbaum--Quillen conjecture involves a filtration on $\text{SH}(F)$ now referred to as the effective slice filtration. This filtration gives rise to a spectral sequence
\[
\text{E}_1 = \uph^{*,*}(F, \mathbb{Z}) \implies \text{K}_*(F)
\]
which computes algebraic K-theory from integral motivic cohomology. 

Both the dual motivic Steenrod algebra and the effective slice spectral sequence are analogues of topological constructions in the stable homotopy category of spectra. The topological dual Steenrod algebra $\aeu^\vee_\mathrm{top}$ serves as the base for the mod-$p$ Adams spectral sequence. The effective slice spectral sequence serves as an analogue of the topological Atiyah--Hirzebruch spectral sequence 
\[
\mathrm{E}_1 = \mathrm{H}^*(X; \mathbb{Z}) \implies \mathrm{KU}_*(X)
\]
which computes complex topological K-theory from singular cohomology.

It was soon realized that many other aspects of stable homotopy theory have generalizations to the motivic setting. Occasionally, these generalizations are relatively straightforward. However, more often than not, these generalizations are much more complicated and rich with topological and algebro-geometric structure.

A central goal in homotopy theory is to understand the stable homotopy groups of spheres $\pi_*(\mathbb{S}^{\mathrm{top}})$. Chromatic homotopy theory organizes these groups into $v_n$-periodic families. Studying the $v_n$-periodicity of $\pi_{*}(\mathbb{S}^{\mathrm{top}})$ is a more approachable task than studying the stable homotopy groups at large, and there are many techniques to do so  \cite{MilRavWil77, Mah81,DevHop04,GoeHenMahRez05}. Moreover, the chromatic filtration extends to the entire stable homotopy category, and the study of its layers has continued to be a focal point of research in modern homotopy theory for over 50 years \cite{Rav84,HS98,BHLS23,BarSchStaWei25Rational}. See \cite{BarBea20} for a detailed survey on the many aspects of chromatic homotopy theory.

In this survey, we will study stable motivic homotopy theory from the chromatic perspective. The central object of study in $\text{SH}(F)$ analogous to $\pi_*(\mathbb{S}^{\mathrm{top}})$ are the stable motivic homotopy groups of spheres $\pi_{*,*}^F(\mathbb{S})$. We will investigate how one can organize $\pi_{*,*}^F(\mathbb{S})$ into periodic families. As we will see, this is a significantly more nuanced story than the topological one. Of particular interest will be computational tools which grant us access to the periodic portion of $\pi_{*,*}^F(\mathbb{S})$. We will look at many variants of the motivic Adams spectral sequence and the slice spectral sequence. 

Additionally, we will keep in mind the ``cofiber of $\tau$" or ``synthetic" mindset. It has been observed that, to some extent, the stable motivic homotopy category $\mathrm{SH}(F)$ is completely governed by toological stable homotopy theory, the Adams--Novikov spectral sequence, and the absolute Galois group of $F$ \cite{HKO-convergencemASS,GIKR22,Pst23,bachmannburklundxu,BHS26}. This philosophy has clarified how to interpret chromatic structures in $\mathrm{SH}(F)$ and given powerful computational tools for computing $\pi_{*}(\mathbb{S}^{\mathrm{top}})$.



\subsection*{Conventions}
We assume that the reader is familiar with the basics of stable homotopy theory, has seen some stable motivic homotopy theory, and that the reader has encountered some aspects of chromatic homotopy theory. We will recall details as is necessary for exposition.

The phrase ``in topology" or ``topologically" will refer to classical stable homotopy theory, i.e. working with the $\infty$-category of spectra $\mathrm{Sp}$ \cite{LurHA}. We will denote the motivic sphere spectrum as $\mathbb{S} \in \mathrm{SH}(F)$ and the topological sphere spectrum as $\mathbb{S}^{\mathrm{top}} \in \mathrm{Sp}$, and other uses of ``top" as a subscript or superscript will refer to working in $\mathrm{Sp}$.

We use bracket notation to denote suspension. For spectra, we let $\mathrm{X}[n] = \Sigma^n \mathrm{X}$, and for motivic spectra we let $\mathrm{X}[s,w] = \Sigma^{s,w}\mathrm{X}$.

\subsection*{Organization}
In \Cref{section:global observations}, we highlight computations of the stable motivic homotopy groups of spheres via the slice spectral sequence. These computations are more geared towards the Milnor--Witt stems at large and are less related to periodicity. In \Cref{section:mass}, we introduce the motivic Adams spectral sequence and synthetic spectra. These tools are the language in which we will study periodicity. In \Cref{section:periodic}, we shift perspectives and look at how topological $v_n$-periodicity manifests in motivic homotopy theory, particularly focusing on rationalization, $v_1$-periodicity, and $v_2$-periodicity. In \Cref{section:exotic}, we study exotic periodic phenomena. We conclude with \Cref{section:future}, where we offer a multitude of future directions of research and open problems.

\subsection*{Acknowledgments}
The author would like to thank Kyle Ormsby and John Palmieri for first teaching him stable, chromatic, and motivic homotopy theory. He would also like to thank Robert Burklund and Piotr Pstragowski for helpful discussions, and Kyle Ormsby and J.D. Quigley for comments and suggestions on an early draft.


\section{Global observations}
\label{section:global observations}

As a graded ring, the topological stable homotopy groups of spheres $\pi_*(\mathbb{S}^{\mathrm{top}})$ are quite rigid. Serre's finiteness theorem implies that $\pi_k(\mathbb{S}^{\mathrm{top}})$ is a finite abelian group for $k>0$ \cite{Serre1953}. As $\pi_0(\mathbb{S}^{\mathrm{top}}) \cong \mathbb{Z}$ and the negative homotopy groups vanish, this allows one to reduce to working at a prime $p$; any integral result may then be obtained by an arithmetic fracture square. In terms of multiplicative structure, Nishida's theorem says that any element of positive degree $\alpha \in \pi_k(\mathbb{S}^{\mathrm{top}})$ must be nilpotent \cite{Nishida73}. 

These structural results dramatically fail in the motivic context. In this section, we investigate large-scale phenomena in $\pi_{*,*}^F(\mathbb{S})$, particularly in comparsion with $\pi_*(\mathbb{S}^{\mathrm{top}})$. Note that in this section, unless otherwise stated, all results are uncompleted.

\subsection{Initial studies}
An early finding of Morel was the isomorphism 
\[
\pi_{0,0}^F(\mathbb{S}) \cong \mathrm{GW}_0(F),
\]
where $\text{GW}_0(F)$ is the \emph{Grothendieck--Witt ring} of symmetric bilinear forms over $F$ \cite{MorelKMW}. Note that $\text{GW}_0(F) = \pi_0(\text{GW}(F))$, where $\text{GW}(F)$ is the (classical, symmetric) hermitian K-theory spectrum of $F$ \cite{Hermitian1}. 

There is more to Morel's finding than just an abstract isomorphism. Key to Voevodsky's solution to the Lichtenbaum--Quillen conjecture was the construction of an $\mathbb{E}_\infty$-ring spectrum $\mathrm{KGL} \in \mathrm{SH}(F)$ called the \emph{algebraic \textup{K}-theory spectrum}. This spectrum represents algebraic $\mathrm{K}$-theory in that there is an isomorphism
\[
\pi_{s,w}^F(\mathrm{KGL}) \cong \mathrm{K}_{s-2w}(F).\footnote{Technically, $\mathrm{KGL}$ only represents Weibel's homotopy algebraic K-theory. However, for regular schemes this agrees with algebraic K-theory \cite{Weibel-KBook}.}
\]
In a similar way, there is an $\mathbb{E}_\infty$-ring spectrum $\mathrm{KQ} \in \mathrm{SH}(F)$ called the \emph{hermitian \textup{K}-theory spectrum} \cite{CalHarNar25}. This spectrum represents hermitian $\mathrm{K}$-theory in that there is an isomorphism
\[
\pi_{s,w}^F(\mathrm{KQ}) \cong \mathrm{GW}_{s-2w}(F)
\]
One may explain Morel's isomorphism as arising from applying $\pi^F_{0,0}(-)$ to the unit map $\mathbb{S} \to \mathrm{KQ}$. In this language, there is more to say.

\begin{thm}[{\cite{MorelKMW}}]
    Let $F$ be any field. The unit map $\mathbb{S} \to \mathrm{KQ}$ induces isomorphisms
    \[
    \pi_{s,w}^F(\mathbb{S}) \cong \left\{\begin{array}{ll}
        \mathrm{K}^{\mathrm{MW}}_{-s}(F) & s=w \\
        0 & s<w.
    \end{array} \right.
    \]
    where $\mathrm{K}^{\mathrm{MW}}_{*}(F)$ is the Milnor-Witt $\mathrm{K}$-theory of $F$.
\end{thm}

Let $\Pi_k^F(\mathbb{S}) = \bigoplus_{n \in \mathbb{Z}}\pi_{n+k, n}^F(\mathbb{S})$ be the $k^{\mathrm{th}}$ \emph{Milnor--Witt stem}. Then Morel's theorem says that there is a grading-reversed isomorphism $\Pi_k^F(\mathbb{S})=0$ for $k<0$ and $\Pi_0^F(\mathbb{S}) \cong \mathrm{K}^{\mathrm{MW}}_*(F)$. 

Milnor--Witt K-theory is an invariant defined via generators and relations in terms of the arithmetic of the field. We recall only a few computations in the context of Morel's theorem, and list the associated values for particular fields in \Cref{table: some values of milnor witt stem 0}. See \cite{Lam05_intro_quadraticforms, Mor12} for more details on Milnor--Witt $\mathrm{K}$-theory. 

\begin{itemize}
        \item As we have seen already,  $\pi_{0,0}^F(\mathbb{S})\cong \mathrm{K}^{\mathrm{MW}}_0(F) \cong \mathrm{GW}_0(F)$, the Grothendieck--Witt ring of quadratic forms over $F$.
        \item We have $\pi_{-1,-1}^F(\mathbb{S}) \cong \mathrm{K}^{\mathrm{MW}}_1(F) \supseteq F^\times$, the units in $F$. 
        \item We have $\pi_{1,1}^F(\mathbb{S}) \cong \mathrm{K}^{\mathrm{MW}}_{-1}(F) \cong \mathrm{W}(F)$, the \emph{Witt ring} of $F$.\footnote{There is an isomorphism $\mathrm{W}(F) \cong \mathrm{GW}_0(F)/(\mathsf{h})$, where $\mathsf{h}$ is the class of the hyperbolic plane over $F$.}
\end{itemize}

\begin{table}[H]
    \centering
    \setlength{\tabcolsep}{0.5em} 
    {\renewcommand{\arraystretch}{1.2}
    \begin{tabular}{|l||l|l|l|}
        \hline
        Field $F$ & $\pi_{-1, -1}^F(\mathbb{S})$ & $\pi_{0,0}^F(\mathbb{S})$ & $\pi_{1,1}^F(\mathbb{S})$ \\
        \hline 
        \hline
        $F=\mathbb{C}$ & contains $\mathbb{C}^{\times}$ \tablefootnote{In particular, this group is never detected upon $p$-completion, hence is invisible to the mod-$p$ motivic Adams spectral sequence.}
        & $\mathbb{Z}$ & $\mathbb{Z}/2$\\
        $F=\mathbb{R}$ & contains $\mathbb{R}^\times$ & $\mathbb{Z} \oplus \mathbb{Z}$\tablefootnote{There is a ring isomorphism $
            \pi_{0,0}^\mathbb{R}(\mathbb{S}) \cong \pi^{C_2}_0(\mathbb{S}_{C_2}) \cong \mathrm{A}(C_2),
            $            where $\pi_\star^{C_2}(\mathbb{S}_{C_2})$ are the $C_2$-equivariant stable stems and $\mathrm{A}(C_2)$ is the Burnside ring.}
        & $\mathbb{Z}$ \\
        $F=\mathbb{Q}_p$, $p \neq 2$ & contains $\mathbb{Q}_p^\times$ & $\mathbb{Z} \oplus \mathbb{Z}/4$ & $\mathbb{Z}/8$ \\
        $F=\mathbb{F}_p$, $p \neq 2$ &  contains $\mathbb{F}_p^\times$ & $\mathbb{Z} \oplus \mathbb{Z}/2$ & $\mathbb{Z}/4$ \\
        \hline
    \end{tabular}}
    \caption{Some particular (additive) values of $\pi_{n,n}^F(\mathbb{S})$.}      
    \label{table: some values of milnor witt stem 0}
\end{table}

Nishida's nilpotence theorem also fails in motivic homotopy theory. For example, there is a Hopf map $\eta \in \pi_{1,1}^F(\mathbb{S})$ coming from the stabilization of the map
\begin{gather*}
S^{3,2} \simeq \mathbb{A}^2 \setminus \{0\} \to  \mathbb{P}^1 \simeq S^{2,1}, \\
(x,y) \mapsto [x:y].
\end{gather*}
This is a motivic analogue of the topological first Hopf map $\eta^{\mathrm{top}} \in \pi_1(\mathbb{S}^\mathrm{top})$, and it was shown by Morel to be non-nilpotent. We will discuss $\eta$-periodicity in further detail in a later section.

\subsection{Slice spectral sequences}
\label{subsection: sliceSS}
A natural next step after Morel's result are the Milnor--Witt stems $\Pi_n^F(\mathbb{S})$ for $n > 0$. To discuss work in this direction, we now provide more details on the slice spectral sequence. 

Let $\mathrm{SH}^{\mathrm{eff}}(F) \subseteq \mathrm{SH}(F)$ denote the full subcategory generated under colimits by objects of the form $\Sigma^{\infty}_{\mathbb{P}^1}X_+[n,0]$ for $X \in \mathrm{Sm}_F$ a smooth scheme over $F$ and $n \in \mathbb{Z}$. We will call this the category of \emph{effective motivic spectra}. The adjoint functor theorem implies that the inclusion of effective motivic spectra admits a right adjoint $r_0$ fitting into a diagram
\[
i_0:\mathrm{SH}^{\mathrm{eff}}(F) \rightleftarrows \mathrm{\mathrm{SH}}(F):r_0.
\]
Let $f_0: = i_0 \circ r_0:\mathrm{SH}(F) \to \mathrm{SH}(F)$ denote the composite. By iterating suspensions, we can define functors $f_q = \Sigma^{q,q}\circ f_0\circ \Sigma^{-q, -q}:\mathrm{SH}(F) \to \mathrm{SH}(F)$ for any $n \in \mathbb{Z}$. For any $\mathrm{E} \in \mathrm{SH}(F)$, we will refer to $f_q(\mathrm{E})$ as the $q^{th}$ \emph{effective cover} of $\mathrm{E}$. When $q=0$, we will say that $f_0(\mathrm{E})$ is the \emph{effective cover} of $\mathrm{E}$.

The $q^{th}$ effective covers assemble naturally into a tower filtering $\mathrm{E}$:
\[\begin{tikzcd}
	\cdots & {f_{q+1}(\mathrm{E})} & {f_q(\mathrm{E})} & {f_{q-1}(\mathrm{E})} & \cdots & {\mathrm{E}}
	\arrow[from=1-1, to=1-2]
	\arrow[from=1-2, to=1-3]
	\arrow[from=1-3, to=1-4]
	\arrow[from=1-4, to=1-5]
	\arrow[from=1-5, to=1-6]
\end{tikzcd}\]
This is the \emph{effective slice tower}. The cofiber $s_q(\mathrm{E}) = \mathrm{cofib}(f_{q+1}(\mathrm{E}) \to f_q(\mathrm{E})$ is the $q^{th}$ \emph{effective slice} of $\mathrm{E}$. Putting this data together gives a diagram of cofiber sequences
\[\begin{tikzcd}
	\cdots & {f_{q+1}(\mathrm{E})} & {f_q(\mathrm{E})} & {f_{q-1}(\mathrm{E})} & \cdots & {\mathrm{E}} \\
	& {s_{q+1}(\mathrm{E})} & {s_q(\mathrm{E})} & {s_{q-1}(\mathrm{E})}
	\arrow[from=1-1, to=1-2]
	\arrow[from=1-2, to=1-3]
	\arrow[from=1-2, to=2-2]
	\arrow[from=1-3, to=1-4]
	\arrow[from=1-3, to=2-3]
	\arrow[from=1-4, to=1-5]
	\arrow[from=1-4, to=2-4]
	\arrow[from=1-5, to=1-6]
\end{tikzcd}\]
Applying motivic homotopy groups gives rise to the \emph{(effective) slice spectral sequence}. This takes the form
\[
\mathrm{E}_1^{s,q,w} = \pi^F_{s,w}(s_q(\mathrm{E})) \implies \pi^F_{s,w}(\mathrm{E}), \quad d_r:\mathrm{E}_r^{s,q,w} \to \mathrm{E}_r^{s-1,q+r, w}.
\]
We will use $\textbf{SliceSS}^F(\mathrm{E})$ as shorthand. The slice spectral sequence converges conditionally to $\mathrm{E}_\eta^\wedge$, the $\eta$-completion of $\mathrm{E}$.\footnote{If $\mathrm{E}$ is any oriented motivic spectrum (or, more generally, a module over an oriented motivic spectrum), then $\mathrm{E} \simeq \mathrm{E}_\eta^\wedge$ \cite{RSOhopf}. For example, $\mathrm{KGL}$ and $\mathrm{kgl}$ are $\eta$-complete, while $\mathrm{KQ}$ and $\mathrm{kq}$ are not.}
\begin{example}
    \label{ex: slicess for kgl}
    Let $\mathrm{E} = \mathrm{KGL}$ be the algebraic K-theory spectrum. By work of Voevodksy and Levine, we have that $s_q(\mathrm{KGL}) \simeq \mathrm{H}\mathbb{Z}[2q,q]$ \cite{Levine08_homotopyconiveau}. Therefore, the $\textbf{SliceSS}^F(\mathrm{KGL})$ takes the form
    \[
    \mathrm{E}_1^{s,q,w} = \bigoplus_{q \in \mathbb{Z}}\pi_{s,w}^F(\mathrm{H}\mathbb{Z}[2q,q]) \implies \pi_{s,w}(\mathrm{KGL}).
    \]
    This allows one to interpret the effective slice spectral sequence as Voevodsky's solution to the Lichtenbaum--Quillen conjecture: it is a spectral sequence which arises from a filtration on algebraic K-theory whose associated graded is motivic cohomology.

    Let $\mathrm{kgl} = f_0(\mathrm{KGL})$ be the effective cover. The slices of $\mathrm{kgl}$ are the same as for $\mathrm{KGL}$ in nonnegative degrees and vanish in negative degrees, so that the $\textbf{SliceSS}^F(\mathrm{kgl})$ acts as a ``connective" version of the $\textbf{SliceSS}^F(\mathrm{KGL})$.\footnote{It is conjectured by Burklund and Krause that up to completion, the slice filtration on $\mathrm{KGL} \in \mathrm{SH}(F)$ coincides with the Hahn, Raksit, and Wilson's even filtration on $\mathrm{K}(F) \in \mathrm{Sp}$ \cite{HRW}. Work in progress of Pstragowski confirms this in the case that $F$ is a local or global field.}
\end{example}

\begin{example}
    \label{ex:slicess for S}
    Let $\mathrm{E} = \mathbb{S}$ be the motivic sphere spectrum. It was conjectured by Voevodsky \cite{Voevodsky02_openproblems} and later proved by Levine \cite{Levine13_convergenceofvoevodskyslicetower} and Rondigs, Spitzweck, and \O stv\ae r \cite{RSOfirst} that the slices of the motivic sphere are related to the topological Adams--Novikov spectral sequence:
    \[
    s_q(\mathbb{S}) \simeq \bigoplus_{k \geq 0}\mathrm{H}(\mathrm{Ext}^{2q-k, k}_{\mathrm{MU}_*\mathrm{MU}}(\mathbb{S}^{\mathrm{top}}))[2q-k, k].
    \]
    Therefore, the $\textbf{SliceSS}^F(\mathbb{S})$ takes the form
    \[
    \mathrm{E}_1^{s,q,w} = \bigoplus_{\substack{q\in\mathbb{Z}\\k \geq 0}}\pi_{s,w}^F(\mathrm{Ext}^{2q-k, k}_{\mathrm{MU}_*\mathrm{MU}}(\mathbb{S}^{\mathrm{top}})) \implies \pi_{s,w}^F(\mathbb{S}).
    \]
\end{example}

\begin{example}
    Let $\mathrm{E}=\mathrm{KQ}$ be the hermitian K-theory spectrum. The slices of hermitian $\mathrm{K}$-theory were determined by R\"ondigs and \O stv\ae r to be \cite{RO16}
    \[
    s_q(\mathrm{KQ}) \simeq \left\{ \begin{array}{rl}
        \mathrm{H}\mathbb{Z}[2q,q] \oplus \bigoplus_{i < \frac{q}{2}}\mathrm{H}\mathbb{F}_2[2i+q,q] &  q \equiv 0 \, (2)\\
        \bigoplus_{i<\frac{q+1}{2}}\mathrm{H}\mathbb{F}_2[2i+q,q] & q \equiv 1 \, (2).
    \end{array} \right.
    \]
    This describes the $\mathrm{E}_1$-page of the $\textbf{SliceSS}^F(\mathrm{KQ})$.

    Let $\mathrm{kq}$ denote the very effective hermitian K-theory spectrum.\footnote{This is the $0^{th}$ cover in the very effective slice filtration of Spitzweck and \O stv\ae r \cite{SpitzweckOstvaer-twistedKtheory}.} Ananyesvkiy, R\"ondigs, and \O stv\ae r determine that the slices of $\mathrm{kq}$ are the same as for $\mathrm{KQ}$ in nonnegative degrees and vanish in negative degrees, so that $\textbf{SliceSS}^F(\mathrm{kq})$ acts as a ``connective" version of $\textbf{SliceSS}^F(\mathrm{KQ})$ \cite{ARO20}.
\end{example}

\begin{remark}
    By the work of Pelaez and of Guti\'errez, R\"ondigs, Spitzweck, and \O stv\ae r, the slice filtration is multiplicative, in that there are pairings
    \[
    f_p(\mathrm{E}) \otimes f_q(\mathrm{E}) \to f_{p+q}(\mathrm{E}), \quad s_p(\mathrm{E}) \otimes s_q(\mathrm{E}) \to s_{p +q}(\mathrm{E}).
    \]
    As a consequence, if $\mathrm{E}$ is a motivic ring spectrum, then $s_0(\mathrm{E})$ is a motivic ring spectrum and $s_*(\mathrm{E}) = \bigoplus_{q \in  \mathbb{Z}}s_q(\mathrm{E})$ is a graded motivic ring spectrum.

    While suppressing some of the delicate multiplicative structure (see \cite{ARO20}, for example), it is conceptually useful to express the graded slices of $\mathrm{kgl}$ and $\mathrm{kq}$ as follows, borrowing notation from \cite{belmontisaksenkong-v1R}:
        \[
        s_*(\mathrm{kgl}) \simeq \mathrm{H}\mathbb{Z}[v_1], \quad s_*(\mathrm{kq}) \simeq \mathrm{H}\mathbb{Z}[\alpha_1, v_1^2]/(2\alpha_1),
        \]
    where $|v_1| = (2,0,1)$ and $\alpha_1 = (1,1,1)$. Here, a monomial in degree $(s,f,w)$ contributes a summand of $\mathrm{HA}[s,w]$ to the $(\frac{s+f}{2})^{th}$-slice, where $\mathrm{A}$ is either $\mathbb{Z}$ or $\mathbb{F}_2$.\footnote{It is worth noting that $\mathrm{kgl}$ and $\mathrm{kq}$ serve as motivic analogues of the topological spectra $\mathrm{ko}$ and $\mathrm{ku}$, that the $\mathrm{E}_2$-pages of the topological Adams--Novikov spectral sequences for $\mathrm{ku}$ and $\mathrm{ko}$ take the form
    \[
    \mathrm{Ext}^{*,*}_{\mathrm{BP}_*\mathrm{BP}}(\mathrm{ku}) \cong \mathbb{Z}[v_1], \quad \mathrm{Ext}^{*,*}_{\mathrm{BP}_*\mathrm{BP}}(\mathrm{ko}) \cong \mathbb{Z}[\alpha_1, v_1^2]/(2\alpha_1),
    \]
    and that the $\mathrm{BP}$-homology of $\mathrm{ku}$ and $\mathrm{ko}$ are are concentrated in even degrees. This observation is well-explained in $\mathrm{SH}(\mathbb{C})$ via the theory of filtered spectra in forthcoming work of Emming.}
\end{remark}

\begin{remark}
    Although we work primarily over fields in this paper, we also remark that Kolderup, Kylling, R\"ondigs, and \O stv\ae r have made extensive use of the slice spectral sequence in $\mathrm{SH}(\mathcal{O}_F[1/2])$, where $\mathcal{O}_F$ is the ring of integers of a number field $F$ \cite{KylRonOst20_kq_quadforms, KolRonOst25_KQoverZ}. 
\end{remark}

\subsection{The first and second Milnor--Witt stems}
The determination of the first Milnor--Witt stem $\Pi_1^F(\mathbb{S})$ is due to Rondigs, Spitzweck, and \O stv\ae r. 

\begin{thm}[\cite{RSOfirst,RSOsecond}]
    For each $n \in \mathbb{Z}$, there is an exact sequence
    \[
    0 \to \mathrm{K}^\mathrm{M}_{2-n}(F)/24 \to \pi_{n+1, n}^F(\mathbb{S}) \to \pi_{n+1, n}^F(\mathrm{kq})\to 0
    \]
    induced by the unit map $\ss \to \mathrm{kq}$.
\end{thm}

For $n=0$, the short exact sequence
\[
    0 \to \mathrm{K}^{M}_{2}(F)/24 \to \pi_{1,0}^F(\mathbb{S}) \to \pi_{1,0}^F(\mathrm{kq}) \to 0
\]
was first conjectured by Morel.

The proof of this result uses an analysis of $\textbf{SliceSS}(\mathbb{S})$ and  $\textbf{SliceSS}(\mathrm{kq})$ and has an interesting relationship (due to the nature of $s_*(\mathbb{S})$) with the topological Adams--Novikov spectral sequence. In subsequent work, the authors build upon their techniques to deduce the second Milnor--Witt stem.

\begin{thm}[\cite{RSOsecond}]
    For each $n \in \mathbb{Z}$, there is an exact sequence
    \[
    0 \to \uph^{1-n, 2-n}(F; \mathbb{Z})/24 \oplus \mathrm{K}^{\mathrm{M}}_{4-n}(F)/2 \to \pi_{n+2, n}^F(\mathbb{S}) \to \pi_{n+2, 2}^F(\upkq)
    \]
    induced by the unit map $\ss \to \upkq$.\footnote{The results of R\"ondigs, Spitzweck, and \O stv\ae r on $\Pi_1^F(\mathbb{S})$ and $\Pi_2^F(\mathbb{S})$ are actually true at the level of homotopy sheaves. These are a stronger invariant than homotopy groups, and $\pi_{*,*}^F(-)$ arise as their global sections. We will not discuss homotopy sheaves any further in this article.}
\end{thm}

Levine has shown that $\pi_{s,0}^\mathbb{C}(\mathbb{S}) \cong \pi_s(\mathbb{S}^{\mathrm{top}})$ \cite{LevineComparison}, and Wilson and \O stv\ae r have shown via different methods that this isomorphism also holds over algebraically closed fields of positive characteristic after inverting the characteristic \cite{WO-finite}.

\section{Motivic Adams spectral sequences}
\label{section:mass}
An alternative method of computation to the slice spectral sequence, and more inspired by chromatic homotopy theory, is the motivic Adams spectral sequence. If $\upe \in \text{SH}(F)$ is a motivic ring spectrum, then it comes equipped with a unit map $\mathbb{S} \to \upe$. The motivic Adams spectral sequence is the tool for computing descent along this map, turning $\mathrm{E}$-homology into stable motivic homotopy groups. 

In more detail, letting $\overline{\upe} = \text{cofib}(\mathbb{S} \to \upe)$, we may form a diagram
\[\begin{tikzcd}
	{\mathbb{S}} & {\overline{\upe}[-1,0]} & {\overline{\upe}[-2,0]^{\otimes 2}} & \cdots \\
	\upe & {\overline{\upe}[-1,0] \otimes \upe} & {\overline{\upe}[-2,0]^{\otimes 2}\otimes \upe}
	\arrow[from=1-1, to=2-1]
	\arrow[from=1-2, to=1-1]
	\arrow[from=1-2, to=2-2]
	\arrow[from=1-3, to=1-2]
	\arrow[from=1-3, to=2-3]
	\arrow[from=1-4, to=1-3]
\end{tikzcd}\]
extending the unit map by iteratively taking fibers.\footnote{Each ``hook" in the diagram is a cofiber sequence.}
For any $\mathrm{X} \in \text{SH}(F)$, applying $\pi_{*,*}^F(\mathrm{X} \otimes -)$ to this diagram results in the E-\emph{based motivic Adams spectral sequence} computing $\pi_{*,*}^F(\mathrm{X})$. This spectral sequence takes the form
\[
\upe_1^{s,f,w} = \pi_{s+f, w}^F(\upe \otimes \overline{\upe}^{\otimes f} \otimes \mathrm{X}) \implies \pi_{s,w}^F(\mathrm{X}_{\upe}^\wedge), \quad d_r:\upe_r^{s,f,w} \to \upe_r^{s-1, f+r, w}.
\]
Here $\mathrm{X}_{\mathrm{E}}^\wedge$ denotes the $\mathrm{E}$-nilpotent completion of $\mathrm{X}$. As in the topological case, some care must be taken on convergence, and only under certain conditions may identify the $\mathrm{E}_2$-page as the cohomology of the Hopf algebroid $(\mathrm{E}_{*,*}(\mathrm{X}), \mathrm{E}_{*,*}(\mathrm{E}\otimes \mathrm{X}))$. 

We now specialize to a few particular cases.

\subsection{$\text{H}\mathbb{F}_2$-based motivic Adams spectral sequence}
The most well-studied motivic Adams spectral sequence is the $\mathrm{H}\mathbb{F}_2$-based motivic Adams spectral sequence. This was first studied by Morel in \cite{Morel99_mAdamsSS}. We will refer to this as \emph{the motivic Adams spectral sequence} for $\mathrm{X}$ and denote it by $\textbf{mASS}^F(\mathrm{X})$. For more details on its construction and convergence, see \cite{DImASS,HKO-convergencemASS, HKORemarks, WO-finite, KylWil19_strongconvergencemASS,Man24}.

The motivic spectrum $\mathrm{H}\mathbb{F}_2$ is flat, and the $\upe_2$-page of the $\textbf{mASS}^F(\mathrm{X})$ takes the form
\[
\upe_2^{s,f,w} = \text{Ext}_{\aeudual}^{s,f,w}(\mathbb{M}_2^F, \uph_{*,*}(\mathrm{X})) \implies \pi_{s,w}^F(\mathrm{X}_{2,\eta}^\wedge).\footnote{Unlike the topological Adams spectral sequence which converges to the $2$-completion, the motivic Adams spectral sequence converges to the $(2,\eta)$-completion. However,by \cite{HKO-convergencemASS}, if $F$ is of finite virtual 2-cohomological dimension, then $\mathrm{X}_{2,\eta}^\wedge \simeq \mathrm{X}_2^\wedge$. }
\]
Here, $\aeu^\vee = \pi_{*,*}^F(\mathrm{H}\mathbb{F}_2 \otimes \mathrm{H}\mathbb{F}_2)$ is the \emph{dual motivic Steenrod algebra}, $\mathbb{M}_2^F = \pi_{*,*}^F(\mathrm{H}\mathbb{F}_2)$ is the mod-2 \emph{motivic homology of a point}, and $\mathrm{H}_{*,*}(\mathrm{X})$ is the mod-2 motivic homology of $\mathrm{X}$. We will omit the $(2,\eta)$-completion from our notation with the assumption that henceforth we are implicitly working in the $(2,\eta)$-completed category $\mathrm{SH}(F)_{2,\eta}^\wedge$. We will often denote the $\mathrm{E}_2$-page of the $\textbf{mASS}^F(\mathrm{X})$ by $\mathrm{Ext}_{\euscr{A}^\vee}^F(\mathrm{X})$.

By Voevodsky's solution to the Milnor conjecture, there is an isomorphism
\[
\pi_{*,*}^F(\mathrm{H}\mathbb{F}_2) \cong (\mathrm{K}^{\mathrm{M}}_*(F)/2)[\tau],
\]
where the \emph{mod-2 Milnor \textup{K}-theory of $F$} $\mathrm{K}^\mathrm{M}_*(F)/2$ is concentrated in the $(n,n)$-degree and $|\tau|=(0, -1)$. We list some values of $\mathbb{M}_2^F$ for various fields in \Cref{table: some values of M2}; see \cite{Lam05_intro_quadraticforms, Weibel-KBook} for more details.

\begin{table}[H]
    \centering
    \begin{tabular}{|c||c|c|c|c|}
        \hline
        &&&&\\[-1em]
        $F$ & $\mathbb{C}$ & $\mathbb{R}$ & $\mathbb{F}_p$, $p\equiv 1 \, (4)$& $\mathbb{F}_p$, $p \equiv 3 \, (4)$  \\
        \hline 
        \hline 
        &&&&\\[-1em]
        $\mathbb{M}_2^F$ & $\mathbb{F}_2[\tau]$&$\mathbb{F}_2[\rho, \tau]$ &$\dfrac{\mathbb{F}_2[u,\tau]}{(u^2)}$ & $\dfrac{\mathbb{F}_2[\rho, \tau]}{(\rho^2)}$\\
        \hline
    \end{tabular}
    
    \vspace{.25cm}
    
    \begin{tabular}{|c||c|c|c|}
        \hline
        &&&\\[-1em]
        $F$ &  $\mathbb{Q}_p$, $p \equiv 1 \, (4)$& $\mathbb{Q}_p$, $p \equiv 3 \, (4)$ & $\mathbb{Q}_2$ \\
        \hline 
        \hline
        &&&\\[-1em]
        $\mathbb{M}_2^F$ & $\dfrac{\mathbb{F}_2[\pi, u, \tau]}{(\pi^2, u^2)}$ & $\dfrac{\mathbb{F}_2[\pi, \rho, \tau]}{(\rho^2, \rho\pi+\pi^2)}$&$\dfrac{\mathbb{F}_2[u,\pi, \rho,\tau]}{(u^2, \pi^2, \rho^3, u\rho, \pi\rho, \rho^2+u\pi)}$\\
        \hline
    \end{tabular}
    \caption{Some values of $\mathbb{M}_2^F$. The degrees of the polynomial generators are $|\rho| = |u| = |\pi|=(-1,-1)$.}     
    \label{table: some values of M2}
\end{table}

By work of Voevodsky \cite{Voemotiviccohomology} and Hoyois, Kelly, and \O stv\ae r \cite{HKO17}, over a field $F$ the dual Steenrod algebra takes the form
\[
\aeu^\vee \cong \mathbb{M}_2^F[\xi_1, \xi_2, \dots, \tau_0, \tau_1, \dots]/(\tau_i^2 +\tau\xi_{i+1}+ \rho\tau_{i+1} + \rho\tau_0\xi_{i+1}),
\]
where $|\xi_{i}| = (2^{i+1}-2, 2^i-1)$ and $|\tau_i| = (2^{i+1}-1, 2^i-1).$ If we need to specify the base field, we will write $\aeu^\vee_F$. Notice that upon inverting $\tau$ and setting $\rho=0$, there is an isomorphism with the topological dual Steenrod algebra
\[
\tau^{-1}\aeu^\vee_{(\rho=0)} \cong \aeu^\vee_{\mathrm{top}}
\]
which sends $\tau_i$ to $\xi_{i+1}^{\mathrm{top}}$ and $\xi_i$ to $(\xi_i^{\mathrm{top}})^2$. As well will see, the element $\tau$ plays a pivotal role in computation in stable homotopy theory.

The pair $(\mathbb{M}_2^F, \aeu^\vee)$ has the structure of a Hopf algebroid, as opposed to how the topological dual Steenrod algebra $\aeu^\vee_{\mathrm{top}}$ has a simpler structure of a Hopf algebra. We will not make direct use of this fact here; see \cite[Appendix A]{Rav86} for more details on Hopf algebroids. 

In contrast with the slice spectral sequence, it is typically much more challenging to compute with the $\textbf{mASS}^F(\mathrm{X})$ without choosing a particular base field $F$. 
We next review some particular base fields and relevant computations with the $\textbf{mASS}^F(\mathbb{S})$.

\subsubsection{$F=\mathbb{C}$}\footnote{By the motivic Lefschetz principle of Balderrama, Culver, and Quigley, there is an equivalence (at least at the level of cellular objects) $\mathrm{SH}(F) \simeq \mathrm{SH}(\mathbb{C})$ for any algebraic closed field $F$ \cite[Proposition 5.2.1]{BalCulQui25}. This equivalence respects the motivic Adams spectral sequences, and so all results in this section translate verbatim to any algebraically closed field.}
\label{subsubsection:motivic adams spectral sequence over C}
The $\textbf{mASS}^\mathbb{C}(\mathbb{S})$ takes the form
\[
\mathrm{E}_2^{s,f,w}=\mathrm{Ext}^\mathbb{C}_{\euscr{A}^\vee}(\mathbb{S}) =\mathrm{Ext}_{\euscr{A}_{\mathbb{C}}^\vee}^{s,f,w}(\mathbb{F}_2[\tau], \mathbb{F}_2[\tau]) \implies \pi_{s,w}^\mathbb{C}(\mathbb{S}).
\]
Key to beginning computations with this spectral sequence is the following result of Dugger ans Isaksen, which is further suggestion that $\tau$-locality relates motivic homotopy theory to topological stable homotopy theory.

\begin{thm}[\cite{DImASS}]
    There is an isomorphism of rings
    \[
    \mathrm{Ext}^{\mathbb{C}}_{\euscr{A}^\vee}(\mathbb{S}) \otimes_{\mathbb{M}_2^\mathbb{C}} \mathbb{M}_2^\mathbb{C}[\tau^{-1}] \cong \mathrm{Ext}_{\euscr{A}^\vee_{\mathrm{top}}}(\mathbb{S}^{\mathrm{top}}) \otimes_{\mathbb{F}_2}\mathbb{F}_2[\tau^{\pm 1}].
    \]
\end{thm}

In \cite{Isaksen19}, Isaksen shows that this isomorphism is one of spectral sequences rather than just $\mathrm{E}_2$-pages. This establishes a compatibility between differentials and hidden extensions in the $\textbf{mASS}^\mathbb{C}(\mathbb{S})$ and the $\textbf{ASS}(\mathbb{S}^{\mathrm{top}})$ that informs many motivic computations from topological ones.

Another relationship between the topological and $\mathbb{C}$-motivic Adams spectral sequences is the following. There is a subalgebra $\euscr{B}^\vee = \mathbb{F}_2[\xi_1, \xi_2, \dots] \subset \aeu^\vee$ and an obvious isomorphism
\[
\aeu^\vee_{\mathrm{top}} \xrightarrow{\cong} \euscr{B}^\vee
\]
where $\xi_i^{\mathrm{top}}$ is sent to $\xi_i$. Notice that $\mathrm{Ext}_{\euscr{B}^\vee}^{*,*,*}(\mathbb{F}_2[\tau], \mathbb{F}_2[\tau])$ is naturally a subalgebra of $\mathrm{Ext}^\mathbb{C}_{\aeu^\vee}(\mathbb{S})$. This is the algebra of elements of Chow degree 0 of $\mathrm{Ext}^\mathbb{C}_{\aeu^\vee}(\mathbb{S})$, where the Chow degree of a class $x$ of degree $(s,f,w)$ is $s+f-2w$. The isomorphism above upgrades to one of Ext groups which respects all higher order structure \cite{Isaksen19}. Hence, there is a subalgebra of $\mathrm{Ext}^{\mathbb{C}}_{\aeu^\vee}(\mathbb{S})$ which is isomorphic to $\mathrm{Ext}_{\euscr{A}^\vee_{\mathrm{top}}}(\mathbb{S}^{\mathrm{top}})$, and the motivic Adams spectral sequence, when restricted to this subalgebra, is isomorphic to the topological Adams spectral sequence.

The isomorphism $\aeu^\vee_{\mathrm{top}} \xrightarrow{\cong}\euscr{B}^\vee$ sends an element of degree $n$ to one of bidegree $(2n, n)$. As a result, at the level of Ext, this isomorphism ``shears" $\mathrm{Ext}_{\euscr{A}^\vee_{\mathrm{top}}}(\mathbb{S}^{\mathrm{top}})$ onto the Chow degree 0 subalgebra of $\mathrm{Ext}^\mathbb{C}_{\euscr{A}^\vee}(\mathbb{S})$.\footnote{For instance, this isomorphism sends $h_j^{\mathrm{top}}$ to $h_{j+1}$.} The precise relationship is that for $s+f-2w=0$, there is an isomorphism \cite{IsaKonLiYuaZhu25}
\[
\mathrm{Ext}^{\frac{s-f}{2}, f}_{\aeu^\vee_{\mathrm{top}}}(\mathbb{F}_2, \mathbb{F}_2) \cong \mathrm{Ext}^{s,f,w}_{\aeu^\vee_\mathbb{C}}(\mathbb{F}_2[\tau], \mathbb{F}_2[\tau]).
\]
Via these methods, as well as the motivic May spectral sequence, $\pi_{*,*}^\mathbb{C}(\mathbb{S})$ was computed through the first 60 stems by Isaksen \cite{Isaksen19}. 

A more powerful technique for computing the $\mathbb{C}$-motivic Adams spectral sequence is the so-called ``cofiber of $\tau$"-method. We will describe this philosophy in further detail in a later section.

\subsubsection{$F=\mathbb{R}$}\footnote{In forthcoming work of Dhankar, Field, Nigam, Quigley, and Yang, a motivic Lefschetz principle is proved, implying that the results in this section translate verbatim to any real closed fields.}
The $\textbf{mASS}^\mathbb{R}(\mathbb{S})$ takes the form
\[
\mathrm{E}_2^{s,f,w} = \mathrm{Ext}^\mathbb{R}_{\euscr{A}^\vee}(\mathbb{S}) = \mathrm{Ext}^{s,f,w}_{\euscr{A}_\mathbb{R}^\vee}(\mathbb{F}_2[\rho,\tau], \mathbb{F}_2[\rho, \tau]) \implies \pi_{s,w}^\mathbb{R}(\mathbb{S}).
\]
The presence of the non-nilpotent $\rho \in \mathbb{M}_2^\mathbb{R}$ makes matters much more complicated than the $\mathbb{C}$-motivic case. However, there is an extremely useful workaround discovered by Hill \cite{Hil11}. Upon filtering the Hopf algebroid $(\mathbb{M}_2^\mathbb{R},\aeu^\vee_\mathbb{R})$ by powers of $\rho$, we obtain an isomorphism of associated gradeds
\[
(\mathrm{gr}^\rho_*\mathbb{M}_2^\mathbb{R}, \mathrm{gr}^\rho_*\aeu^\vee_\mathbb{R}) \cong (\mathbb{M}_2^\mathbb{C}, \aeu^\vee_\mathbb{C}) \otimes_{\mathbb{F}_2} \mathbb{F}_2[\rho]
\]
This filtration lifts to a filtration on the cobar complex computing $\mathrm{Ext}$, hence gives a $\rho$-\emph{Bockstein spectral sequence}
\[
\mathrm{E}_1 = \mathrm{Ext}^{\mathbb{C}}_{\aeu^\vee}(\mathbb{S})[\rho] \implies \mathrm{Ext}^\mathbb{R}_{\aeu^\vee}(\mathbb{S}).
\]
In \cite{DugIsa17}, Dugger and Isaksen use the $\rho$-Bockstein to compute with the $\textbf{mASS}^\mathbb{R}(\mathbb{S})$, where they determine the first four Milnor--Witt stems over $\mathbb{R}$. This technique was later extended by Belmont and Isaksen to compute up to the $11^{th}$ Milnor--Witt stem \cite{BelIsa22-rmotivicstems}.

In a different direction, Balderamma, Culver, and Quigley have fully computed $\mathrm{Ext}^\mathbb{R}_{\aeu^\vee}(\mathbb{S})$ in filtrations $f \leq 3$ by constructing and using a motivic lambda algebra \cite{BalCulQui25}.

There is also a substantial interaction between $\mathbb{R}$-motivic and $C_2$-equivariant homotopy theory that has led to much progress in the computation of the $\textbf{mASS}^\mathbb{R}(\mathbb{S})$ \cite{HelOrm-Galois,DugIsa17_rmotivicC2equivariant,BelGuiIsa21_C2equivariantRmotivic2, GuiIsa24_c2stems}.

\subsubsection{$F=\mathbb{F}_p$ for $p \neq 2$.}
The $\textbf{mASS}^{\mathbb{F}_p}(\mathbb{S})$ takes the form
\[
\mathrm{E}_2^{s,f,w} = \mathrm{Ext}^{\mathbb{F}_p}_{\aeu^\vee}(\mathbb{S}) = \mathrm{Ext}^{s,f,w}_{\euscr{A}_{\mathbb{F}_p}^\vee}\left(\dfrac{\mathbb{F}_2[x, \tau]}{(x^2)}, \dfrac{\mathbb{F}_2[x, \tau]}{(x^2)} \right) \implies \pi_{s,w}^{\mathbb{F}_p}(\mathbb{S}),
\]
where $x=u$ if $p \equiv 1 \, (4)$ and $x = \rho$ if $p \equiv 3 \, (4)$. This Adams spectral sequence was the particular study of Wilson and \O stv\ae r in \cite{WO-finite}. They compute $\pi_{n,0}^{\mathbb{F}_p}(\mathbb{S})$ for $0 \leq n \leq 18$, curiously producing an isomorphism
\[
\pi_{n,0}^{\mathbb{F}_p}(\mathbb{S}) \cong \pi_{n}(\mathbb{S}^{\mathrm{top}}) \oplus \pi_{n+1}(\mathbb{S}^{\mathrm{top}}),
\]
although this is known to not hold in general. 

It is worth noting that there are significantly more differentials in the $\textbf{mASS}^{\mathbb{F}_p}(\mathbb{S})$. For example, while the element $\tau$ lives on the $\mathrm{E}_2$-page, one can show that no power of $\tau$ can live in $\pi_{*,*}^{\mathbb{F}_p}(\mathbb{S})$.\footnote{For example, the map $\mathbb{S} \to \mathrm{H}\mathbb{Z}$ induces a map of motivic Adams spectral sequences which is surjective on $\mathrm{E}_2$-pages in stems 0 and -1. The integral motivic cohomology groups of $\mathbb{F}_p$ do not allow for any $\tau$-divisibility in stem 0 on the $\mathrm{E}_\infty$-page of the $\textbf{mASS}^{\mathbb{F}_p}(\uph\mathbb{Z})$, and their particular form uniquely determines the differentials which kill all powers of $\tau$ \cite{Soule79,WO-finite}. By naturality, these differentials pull back to the $\textbf{mASS}^{\mathbb{F}_p}(\mathbb{S}).$} By the Liebniz rule, there must infinitely many pages which have differentials originating from the stem 0 and filtration 0 coordinate to kill these powers of $\tau$. 

\begin{remark}
    There is a common thread found in computation over any base field of ``reducing to a simpler Ext group". One way that this manifests in $\mathbb{R}$-motivic computations in that there is an isomorphism due to Belmont and Isaksen $
\mathrm{Ext}^\mathbb{R}_{\euscr{A}^\vee}(\mathbb{S}/\rho) \cong \mathrm{Ext}^\mathbb{C}_{\aeu^\vee}(\mathbb{S}),$
where $\mathbb{S}/\rho = \mathrm{cofib}(\mathbb{S}[-1,-1] \xrightarrow{\rho}\ss).$ \cite{BelIsa22-rmotivicstems} More generally, Behrens and Shah have shown that there is an equivalence of categories $\mathrm{SH}(\mathbb{C}) \simeq \mathrm{Mod}_{\mathbb{S}/\rho}(\mathrm{SH}(\mathbb{R}))$ \cite{BS20}.

    Another useful insight, due to Balderrama, Culver, and Quigley \cite[Remark 2.2.8]{BalCulQui25}, is that for any field $F$ there is an isomorphism
    \[
    \mathrm{Ext}_{\aeu^\vee}^F(\mathbb{S}) \cong \mathrm{Ext}^\mathbb{R}_{\aeu^\vee}(\mathbb{M}_2^\mathbb{R}, \mathbb{M}_2^F).
    \]
    This philosophy is, in spirit, a generalization of the $\rho$-Bockstein spectral sequence. In this way, it is possible to compute $\mathrm{Ext}^{F}_{\aeu^\vee}(\mathbb{S})$ by knowing $\mathrm{Ext}^{\mathbb{R}}(\mathbb{S})$ and the $\rho$-action on $\mathbb{M}_2^F$. The previous authors use this idea to make computations with the $\mathrm{Ext}^{\mathbb{Q}_p}_{\aeu^\vee}(\mathbb{S})$. These ideas extend to other spectra than the sphere, and are a key input to forthcoming work of the author, Petersen, and Tatum \cite{MorPetTat}.
\end{remark}

\subsection{Motivic Adams--Novikov spectral sequence}
Another important motivic Adams spectral sequence is the motivic Adams--Novikov spectral sequence. Let $\mathrm{MGL}$ be the \emph{algebraic cobordism spectrum} \cite{LM07}. By work of Hu and Kriz, and of Vezzosi, one can emulate Quillen's idempotent construction to split this spectrum into a wedge of suspensions of the spectrum $\mathrm{BPGL}$, the \emph{motivic Brown--Peterson spectrum} \cite{Hu-Kriz-remarks,Vez01_motivicBP}. The $\mathrm{BPGL}$-based motivic Adams spectral sequence is what we will refer to as the \emph{motivic Adams--Novikov spectral sequence} and denote by $\textbf{mANSS}^F(\mathrm{X})$. This spectral sequence has primarily been studied over $\mathbb{C}$. For more details on its construction and convergence, see \cite{DImASS, HKORemarks}.

The motivic spectrum $\mathrm{BPGL}$ is flat, and the $\mathrm{E}_2$-page of the $\textbf{mANSS}^F(\mathrm{X})$ takes the form
\[
\mathrm{E}_2^{s,f,w} = \mathrm{Ext}_{\mathrm{BPGL}_{*,*}\mathrm{BPGL}}^{s,f,w}(\mathrm{BPGL}_{*,*}, \mathrm{BPGL}_{*,*}(\mathrm{X})) \implies \pi_{s,w}^F(\mathrm{X}).
\]
We will often denote the $\mathrm{E}_2$-page of the $\textbf{mANSS}^F(X)$ by $\mathrm{Ext}_{\mathrm{BPGL}}^F(\mathrm{X})$.

For $F=\mathbb{C}$, Hu, Kriz, and Ormsby showed that there is an isomorphism of rings \cite{HKORemarks}
\[\mathrm{Ext}_{\mathrm{BPGL}}^{\mathbb{C}}(\mathbb{S}) \cong \mathrm{Ext}_{\mathrm{BP}_*\mathrm{BP}}(\mathbb{S}^{\mathrm{top}})\otimes \mathbb{Z}[\tau].\footnote{One must appropriately regrade the topological Adams--Novikov $\mathrm{E}_2$-page into Chow degree 0 for the grading on this isomorphism to work.}\]
Notice that, as opposed to the findings of Dugger and Isaksen with the $\textbf{mASS}^\mathbb{C}(\mathbb{S})$, one does not need to invert $\tau$ on the motivic Adams--Novikov $\mathrm{E}_2$-page for this equivalence to hold. As localization is quite destructive, this implies that the relationship between the $\textbf{mANSS}^\mathbb{C}(\mathbb{S})$ and the $\textbf{ANSS}^\mathbb{C}(\mathbb{S}^{\mathrm{top}})$ is much closer than that of the $\textbf{mASS}^\mathbb{C}(\mathbb{S})$ and the $\textbf{ASS}(\mathbb{S}^{\mathrm{top}})$.

\subsection{Synthetic spectra}

The connection between the $\mathbb{C}$-motivic Adams--Novikov spectral sequences and the topological Adams--Novikov spectral sequences suggests that there is an intricate relationship between the stable motivic homotopy category $\mathrm{SH}(\mathbb{C})$ and the topological category of spectra $\mathrm{Sp}$, and that this relationship is governed in some way by $\tau \in \pi_{0, -1}^{\mathbb{C}}(\mathbb{S})$. Indeed, there is much to say. For the rest of this section, we restrict our attention to the cellular subcategory $\mathrm{SH}(F)^{\mathrm{cell}} \subseteq \mathrm{SH}(F)$ unless otherwise stated \cite{DI05}. This is the smallest subcategory of $\mathrm{SH}(F)$ generated under colimits by $\mathbb{S}[s,w]$ for $s,w \in \mathbb{Z}$.

There is a Betti realization functor $\text{Be}^\mathbb{C}:\text{SH}(\mathbb{C}) \to \text{Sp}$
which is a stable incarnation of the functor sending a scheme $X \in \mathrm{Sm}_\mathbb{C}$ to its space $X(\mathbb{C})$ of $\mathbb{C}$-points.\footnote{If $F$ embeds into $\mathbb{C}$, then there is an extension of $\mathrm{Be}^\mathbb{C}$ to $\mathrm{SH}(F)$ taking the same  values. By work of Bachmann and \O stv\ae r one can also make use of this functor when $F$ has positive characteristic \cite{BO22}.} We list some values of the Betti realization in \Cref{table: betti realization}. 

\begin{table}[H]
    \centering
    \setlength{\tabcolsep}{0.5em} 
    {\renewcommand{\arraystretch}{1.2}
    \begin{tabular}{|c||c|}
        \hline
        $\mathrm{E} \in \mathrm{SH}(\mathbb{C})$ & $\mathrm{Be}^\mathbb{C}(\mathrm{E}) \in \mathrm{Sp}$\\
        \hline 
        \hline
        $\mathbb{S}$ & $\mathbb{S}^{\mathrm{top}}$\\
        $\mathrm{HA}$ & $\mathrm{HA}^{\mathrm{top}}$\\
        $\mathrm{KGL}, \mathrm{kgl}$ & $\mathrm{KU}, \mathrm{ku}$\\
        $\mathrm{KQ}, \mathrm{kq}$ & $\mathrm{KO}, \mathrm{ko}$\\
        $\mathrm{MGL}$ & $\mathrm{MU}$\\
        $\mathrm{BPGL}$, $\mathrm{BPGL}\langle n \rangle$ & $\mathrm{BP}$, $\mathrm{BP}\langle n \rangle$\\
        \hline
    \end{tabular}}
    \caption{Some values of $\mathrm{Be}^\mathbb{C}(\mathrm{E})$.}      
    \label{table: betti realization}
\end{table}

It was shown by Dugger and Isaksen \cite{DImASS} that Betti realization has the effect of ``inverting $\tau$''; this is the functor which realizes the aforementioned isomorphism $\tau^{-1}\aeu^\vee \cong \aeu^\vee_{\mathrm{top}}$ of dual Steenrod algebras. They further identified this $\tau$-localized category as the topological category of spectra
\[
\tau^{-1}\mathrm{SH}(\mathbb{C}) \simeq \mathrm{Sp}.
\]
In other words, setting ``$\tau=1$" recovers topological stable homotopy theory. 

On the other hand, one may set ``$\tau=0$". Let $\mathbb{S}/\tau = \mathrm{cofib}(\mathbb{S} [0,-1]\xrightarrow{\tau}\ss)$ be the cofiber of $\tau$.
Surprisingly, Gheorghe showed that $\mathbb{S}/\tau \in \mathrm{CAlg}(\mathrm{SH}(\mathbb{C}))$ \cite{Gheorghe18}.\footnote{We say surprising because Moore spectra often have terrible multiplicative properties \cite{Rondigs20_motivicmoorespectra,Burklund22_moorespectra}.} Even more surprisingly, Isaksen showed that there is an isomorphism of graded abelian groups \cite{Isaksen19}
\[
\pi_{*,*}^\mathbb{C}(\mathbb{S}/\tau) \cong \text{Ext}_{\text{BP}_*\text{BP}}^{*,*}(\mathbb{S}^{\mathrm{top}}).
\]
In other words the $\upe_2$-page of the $\textbf{ANSS}(\mathbb{S}^{\mathrm{top}})$ is precisely encoded by the homotopy of the motivic spectrum $\mathbb{S}/\tau$. 

Let us tighten the thread a little further. A useful tool for computing the $\mathrm{E}_2$-page of the $\textbf{ANSS}(\mathbb{S}^{\mathrm{top}})$ is the algebraic--Novikov spectral sequence, which arises from filtering $\mathrm{BP}_*\mathrm{BP}$ by powers of the ideal $I=(p, v_1, v_2, \dots) \subset \pi_*(\mathrm{BP})$.\footnote{The $\mathrm{E}_2$-page of the algebraic Novikov spectral sequence involves computing Ext over the algebra $\euscr{P}^\vee = \mathrm{H}_*(\mathrm{BP}) \cong \mathbb{F}_2[\xi_1^2, \xi_2^2, \dots] \subset \aeu^\vee_{\mathrm{top}}$, which we refer to as the Palmieri subalgebra.} We will refer to this spectral sequence by $\textbf{aNSS}(\mathbb{S}^{\mathrm{top}}).$ By the identification of Isaksen, another spectral sequence computing the $\mathrm{E}_2$-page of the $\textbf{ANSS}(\mathbb{S}^{\mathrm{top}})$ is the $\mathbf{mASS}^\mathbb{C}(\mathbb{S}/\tau)$. 

\begin{thm}[\cite{GheWanXu21}]
    There is an isomorphism of spectral sequences, up to reindexing, between the $\textup{\textbf{mASS}}^\mathbb{C}(\mathbb{S}/\tau)$ and the $\textup{\textbf{aNSS}}(\mathbb{S}^{\mathrm{top}})$.
    \[
    \begin{tikzcd}
	{\mathrm{Ext}_{\aeu^\vee}^{\mathbb{C}}(\mathbb{S}/\tau)} & & {\pi_{*,*}^\mathbb{C}(\mathbb{S}/\tau)} \\
	{\mathrm{Ext}_{\euscr{P}^\vee}(I^*/I^{*+1})} &  &{\mathrm{Ext}_{\mathrm{BP}_*\mathrm{BP}}(\mathbb{S}^{\mathrm{top}})}
	\arrow[Rightarrow,"{\textup{\textbf{mASS}}^\mathbb{C}(\mathbb{S}/\tau)}", from=1-1, to=1-3]
	\arrow[equals,from=1-1, to=2-1]
	\arrow[equals,from=1-3, to=2-3]
	\arrow[Rightarrow,"{\textup{\textbf{aNSS}}(\mathbb{S}^{\mathrm{top}})}"', from=2-1, to=2-3]
    \end{tikzcd}
    \]
\end{thm}

Gheorghe, Wang, and Xu actually showed something much stronger \cite{GheWanXu21}. Just as $\text{SH}(\mathbb{C}) \simeq \text{Mod}_{\ss}(\text{SH}(\mathbb{C}))$, one may consider the category $\text{Mod}_{\ss/\tau}(\text{SH}(\mathbb{C}))$ of modules over $\ss/\tau$.
They show that there is an equivalence of categories
\[
\mathrm{Mod}_{\ss/\tau}(\mathrm{SH}(\mathbb{C})) \simeq \mathrm{Stable}^{\mathrm{ev}}(\mathrm{BP}_*\mathrm{BP}),
\]
where $\mathrm{Stable}^{\mathrm{ev}}(\mathrm{BP}_*\mathrm{BP})$ is the stable derived category of even comodules over $\mathrm{BP}_*\mathrm{BP}$. 

These realizations have lead to dramatic changes in the way that the motivic Adams spectral sequence is used to compute both $\pi_*(\mathbb{S}^{\mathrm{top}})$ and $\pi_{*,*}^\mathbb{C}(\mathbb{S})$. The rough idea is as follows (see also \cite[Section 4.3]{IsaOst20}).
We would like to compute $\pi_{*,*}^\mathbb{C}(\mathbb{S})$ by the $\textbf{mASS}^{\mathbb{C}}(\mathbb{S})$. Begin by computing the $\mathrm{E}_2$-page, which we may do by machine in a large range, and resolve differentials using any previously mentioned technique.
Next, pass to $\mathbb{S}/\tau$. The $\textbf{mASS}^\mathbb{C}(\mathbb{S}/\tau)$ is isomorphic to the $\textbf{aNSS}(\mathbb{S}^{\mathrm{top}})$, which we may also compute by machine in a large range.
Using the cofiber sequence $\mathbb{S} \to \mathbb{S}/\tau \to \mathbb{S}[1,-1]$ and the fact that the motivic Adams spectral sequence is natural, we may pull back and push forward differentials in the $\textbf{mASS}^\mathbb{C}(\mathbb{S}/\tau)$ to the $\textbf{mASS}^\mathbb{C}(\mathbb{S})$. 

The key point is that the $\textbf{aNSS}(\mathbb{S}^{\mathrm{top}})$ is more easily computable by machine as it is algebraic in nature. Moreover, by inverting $\tau$ on the end result, we obtain results about the $\textbf{ASS}(\mathbb{S}^{\mathrm{top}})$ and the $\pi_*(\mathbb{S}^{\mathrm{top}}).$ It is by this technique that the knowledge of $\pi_{*,*}^\mathbb{C}(\mathbb{S})$ has been pushed furthest. Isaksen, Wang, and Xu have use the ``cofiber of $\tau$"-method to compute up to stem 90 in both $\pi_{*,*}^\mathbb{C}(\mathbb{S})$ and $\pi_*(\mathbb{S}^{\mathrm{top}}).$ \cite{IWX23}

One can organize this relationship into the following diagram:
\[
\mathrm{Sp} \xleftarrow{\tau^{-1}}\mathrm{SH}(\mathbb{C}) \xrightarrow{ \otimes \ss/\tau} \mathrm{Stable}^{\mathrm{ev}}(\mathrm{BP}_*\mathrm{BP}).
\]
In this way, one can view the topological stable homotopy category as the generic fiber and the derived category over $\text{BP}_*\text{BP}$ as the special fiber of a family of categories over $\mathbb{A}^1$ which deform $\mathrm{SH}(\mathbb{C})$ and are parameterized by $\tau$.\footnote{Work of Bachmann, Kong, Wang, and Xu have constructed a $t$-structure on $\mathrm{SH}(F)$ for any field $F$, known as the Chow $t$-structure, such that the heart $\mathrm{SH}(F)^{\heartsuit}$ of this $t$-structure is equivalent to the stable derived categoy of $\mathrm{BP}_*\mathrm{BP}$-comodules \cite{BacKonWanXu_Chowtstructure}.} 

This was greatly generalized to the theory of synthetic spectra by Pstragowski \cite{Pst23}. To any homotopy commutative and associative ring spectrum $\text{E} \in \text{Sp}$ which is Adams-type,
Pstragowski constructed the category $\text{Syn}_\text{E}$ of $\text{E}$-synthetic spectra. There is a canonical element $\tau \in \pi_{0, -1}^\mathrm{E}(\mathbb{S}_\mathrm{E})$ in the homotopy of the unit object $\mathbb{S}_\mathrm{E}$. In parallel with the $\mathbb{C}$-motivic story, one can set $\tau=1$ and $\tau=0$, and the resulting span interpolates between spectra and the derived category over the cooperations $\mathrm{E}_*\mathrm{E}$:
\[
\mathrm{Sp} \xleftarrow{\tau^{-1}} \mathrm{Syn}_\mathrm{E} \xrightarrow{\otimes \ss_\mathrm{E}/\tau} \mathrm{Stable}(\mathrm{E}_*\mathrm{E}).
\]
Pstragowski demonstrated that there is an equivalence of categories $\text{SH}(\mathbb{C}) \simeq \text{Syn}_{\text{BP}},$ and the deformation perspectives are also isomorphic.

The category $\mathrm{Syn}_{\mathrm{E}}$ is equipped with a synthetic analogue functor
\[
\nu_\mathrm{E}:\mathrm{Sp} \to \mathrm{Syn}_\mathrm{E}
\] 
which acts as a section to $\tau^{-1}$.\footnote{Since one may identify $\mathrm{Be}^\mathbb{C}$ with $\tau^{-1}$, it follows that $\nu_{\mathrm{BP}}$ is a section of $\mathrm{Be}^{\mathbb{C}}$, so one may reinterpret \Cref{table: betti realization} as topological spectra and their $\mathrm{BP}$-synthetic (or $\mathbb{C}$-motivic) analogues.} The synthetic analogue $\nu_{\mathrm{E}}(\mathrm{X})$ encodes both the homotopy groups $\pi_*(\mathrm{X})$ and the data of the $\mathrm{E}$-$\textbf{ASS}(\mathrm{X})$. There is an isomorphism 
\[
\pi_{*,*}^\mathrm{E}(\nu_\mathrm{E}(\mathrm{X})/\tau) \cong \mathrm{Ext}_{\mathrm{E}_*\mathrm{E}}^{*,*}(\mathrm{X})
\]
between the synthetic homotopy groups of $\nu_\mathrm{E}(\mathrm{X})/\tau$ and the $\mathrm{E}_2$-page of the $\mathrm{E}$-$\textbf{ASS}(\mathrm{X})$. Let $^{\mathrm{syn}}\mathrm{E}_2^{*,*,*}(\nu_\mathrm{E}(\mathrm{X}))$ denote the $\mathrm{E}_2$-page of the $\nu_\mathrm{E}(\mathrm{E})$-based synthetic Adams spectral sequence for $\nu_\mathrm{E}(\mathrm{X})$ and $\mathrm{E}_2^{*,*}(\mathrm{X})$ denote the $\mathrm{E}_2$-page of the $\mathrm{E}$-$\textbf{ASS}(\mathrm{X})$. Then 
\[^{\mathrm{syn}}\mathrm{E}_2^{*,*,*}(\nu_\mathrm{E}(\mathrm{X})) \cong \mathrm{E}_2^{*,*}(\mathrm{X}) \otimes \mathbb{Z}[\tau]\] Moreover, if there is a topological differential $d_r^{\mathrm{top}}(x)=y$ in the $\mathrm{E}$-$\textbf{ASS}(\mathrm{X})$, then there is a differential $d_r^{\mathrm{syn}}(x) = \tau^{r-1}y$ in the synthetic analogue.\footnote{In the case of $\mathrm{E} = \mathrm{BP}$, this recovers the observations of Hu, Kriz, and Ormsby in the $\textbf{mANSS}^\mathbb{C}(\mathbb{S})$ \cite{HKORemarks}.}


    

An analogous topological reconstruction of $\mathrm{SH}(\mathbb{C})$ was done by Gheorghe, Isaksen, Krause, and Ricka using filtered spectra \cite{GIKR22}. They construct a functor $\Gamma_\star:\mathrm{Sp} \to \mathrm{Fil}(\mathrm{Sp})$ which takes a spectrum $\mathrm{X}$ to the  decalage of the Adams--Novikov filtration for $\mathrm{X}$, in that
\[
    \Gamma_w(\mathrm{X}) \simeq \mathrm{Tot}(\tau_{\geq 2w}(\mathrm{X}\otimes \mathrm{MU}^{\otimes \bullet +1})).
\]
The authors then show that there is an equivalence of categories 
\[\mathrm{Mod}_{\Gamma_\star( \mathbb{S}^{\mathrm{top}})}(\mathrm{Fil}(\mathrm{Sp})) \simeq \mathrm{SH}(\mathbb{C}).
\] 
As filtered spectra are naturally intertwined with spectral sequences, they are typically more amenable to computation than synthetic spectra. Throughout the literature, one can find many reinterpretations of $\mathrm{Syn}_\mathrm{E}$ via filtered models \cite{BurXu_hj3,BaeJohMar25_stablecomoduledeformations,CarDavNig25_descentssviasynthetic}. One can think of these filtered models as taking a spectrum $\mathrm{X}$ to the decalage of the $\mathrm{E}$-Adams filtration for $\mathrm{X}$

The synthetic perspective of $\mathrm{SH}(F)$ has seen development beyond the case of $F=\mathbb{C}$. In the case of $F=\mathbb{R}$, Burklund, Hahn, and Senger use the $\text{Gal}(\mathbb{C}/\mathbb{R})=C_2$-action on the $\mathbb{C}$-points of a real motivic spectrum to interpret $\text{SH}(\mathbb{R})$ as a one-parameter deformation of $\text{Sp}^{C_2}$ \cite{BHS26}.\footnote{If $F\subseteq \mathbb{R}$, there is a Real Betti realization functor $\mathrm{Be}^\mathbb{R}:\mathrm{SH}(F) \to \mathrm{Sp}$. This functor is far different than its $\mathbb{C}$-motivic counterpart; see \cite{Bann-realbetti} for a wonderful recollection.} The special fiber is related to the stable derived category of $\mathrm{MU}_{\mathbb{R}\star}\mathrm{MU}_\mathbb{R}$-comodules.

Angelini-Knoll, Behrens, Belmont, and Kong have attempted to replicate these arguments and create an ``artificial'' motivic category which is a deformation of $\mathrm{Sp}^{C_p}$ for $p$ odd \cite{AngKnoBehBelKon_DeformationofBorel}. A hurdle to connecting motivic homotopy theory with $C_p$-equivariant homotopy theory is that, while $\mathrm{Gal}(\mathbb{C}/\mathbb{R})\cong C_2$, for $p$ odd there is no algebraically closed field $F$ with a subfield $L$ such that $\mathrm{Gal}(F/L) \cong C_p$. 

There have also been applications of the ideas of synthetic spectra in the context of equivariant motivic homotopy theory. For $F=\mathbb{C}$ and $\mathrm{A}$ a finite abelian group, Allen and Piessevaux have shown that the category $\mathrm{SH}^{\mathrm{A}}(\mathbb{C})$ of $\mathrm{A}$-equivariant $\mathbb{C}$-motivic spectra admits a deformation interpretation via an equivariant lift of $\tau$ \cite{AllPie25_SynSHCequivariant}. The authors show that the generic fiber may be identified with the category $\mathrm{Sp}^\mathrm{A}$ of global $\mathrm{A}$-equivariant spectra, and the special fiber is related to the stable derived category of $\mathrm{MU}_{\mathrm{A}\star}\mathrm{MU}_\mathrm{A}$-comodules.\footnote{Curiously, Allen and Piessevaux show that the deformation parameter $\tau$ is divisible by various equivariant enhancements $\tau_\alpha$, where $\alpha$ is a character of $\mathrm{A}$. Work in progress of the previously mentioned authors identifies the categories of modules over $\mathbb{S}/\tau_\alpha$.} See also the work of Carrick on filtered $G$-spectra \cite{Carrick25_slicesynthetic}.

In recent work of Bachmann, Burklund, and Xu, where they reconstruct $\mathrm{SH}(F)$ for fields $F$ of small cohomological dimension in terms of their absolute Galois groups \cite{bachmannburklundxu}. As a consequence, they show that for any field $F$ containing all $2^{nd}$ roots of unity that there is an isomorphism
\[
\pi_{*,*}^F(\mathbb{S}) \cong \pi_{*,*}^{\text{BP}}(\mathbb{S}_{\mathrm{BP}}) \otimes \mathrm{K}^{\mathrm{MW}}_*(F),\footnote{Work in progress of Klaus Mattis and Anton Englemann constructs a synthetic reconstruction in the setting of \'etale motivic spectra.}
\]
where $\mathbb{S}_{\mathrm{BP}}=\nu_{\mathrm{BP}}(\mathbb{S}^{\mathrm{top}})$ and $\pi_{*,*}^{\mathrm{BP}}(-)$ denotes homotopy groups in $\mathrm{Syn}_{\mathrm{BP}}$.

Finally, we note that recent work of Tanania, building on observations of Burklund, Hahn, and Senger, has shown a surprising connection between $\mathbb{F}_2$-synthetic spectra and $\mathbb{R}$-motivic homotopy theory \cite{tanania2025_realisotropic}. Tanania shows that there is a ring spectrum $\mathbb{S}^\mathrm{iso} \in \mathrm{CAlg}(\text{SH}(\mathbb{R}))$ called the isotropic sphere spectrum, and that there is an equivalence of categories
\[
\text{Mod}_{\mathbb{S}^{\mathrm{iso}}}(\mathrm{SH}(\mathbb{R})) \simeq \text{Syn}_{\mathbb{F}_2}.
\]
Moreover, he shows that the element $\rho$ lifts to an element of $\pi_{-1, -1}^\mathbb{R}(\mathbb{S}^{\mathrm{iso}})$ and serves as a deformation parameter, in that there is a span
\[
\text{Sp} \xleftarrow{\rho^{-1}}\text{SH}^\text{iso}(\mathbb{R}) \xrightarrow{ \otimes \mathbb{S}^\text{iso}/\rho} \mathrm{Stable}(\euscr{A}^\vee_{\mathrm{top}}).\footnote{It was shown by Bachmann that $\rho^{-1}\mathrm{SH}(\mathbb{R}) \simeq \mathrm{Sp}$ without the presence of isotropic localization \cite{Bachmann-real}.}
\]
Notice that $\euscr{A}^\vee_{\mathrm{top}} = \pi_{*}(\mathrm{H}\mathbb{F}_2^{\mathrm{top}}\otimes \mathrm{H}\mathbb{F}_2^{\mathrm{top}})$ is the topological dual Steenrod algebra, hence the above span is the same one gets for $\mathrm{Syn}_{\mathbb{F}_2}$ with parameter $\tau$.

\begin{remark}
    It is tantalizing to interpret the $\mathbb{F}_2$ appearing in $\mathrm{Syn}_{\mathbb{F}_2}$ in Tanania's work as $C_2 = \mathrm{Gal}(\mathbb{C}/\mathbb{R})$, which suggests that in the artificial motivic categories of Angelini-Knoll, Behrens, Belmont, and Kong, there is a subcategory which is equivalent to $\mathrm{Syn}_{\mathbb{F}_p}$.
\end{remark}

\section{Periodic phenomena}
\label{section:periodic}

We now shift perspectives to study periodicity. Although the ring structure on $\pi_*(\mathbb{S}^{\mathrm{top}})$ is ill-behaved, there are deeper structures which are revealed by chromatic homotopy theory. It was first observed by Adams that for $i=1,2$ and all $k>0$, there are infinite families of elements in the stable homotopy groups which we denote $v_1^{4k}(\eta^i) \in \pi_{8k+i}(\mathbb{S}^{\mathrm{top}})$. 

The notation $v_1^{4k}(-)$ is nonstandard,\footnote{Other names for these families include $P^k\eta^i$ and $\mu_{8k+i}$.} but illustrates a crucial observation of Adams: there is an operator $v_1$ which acts on $\pi_*(\mathbb{S})$ non-nilpotently, generating infinite families of elements. This operator is not itself an element of $\pi_*(\mathbb{S}^{\mathrm{top}})$, and our use of parentheses should indicate that this not an actual product of stable homotopy elements. This operator fits into a family of operators $v_n$ for all $n \geq 0$, and every element in $\pi_*(\mathbb{S}^{\mathrm{top}})$ fits into some $v_n$-periodic family. 

There is a rich connection between the $v_n$-operators and arithmetic geometry, specifically via the moduli stack of formal groups. We will instead restrict our attention primarily to computational techniques and not divert our attention into this area. The interested reader should consult \cite{Rav86,COCTALOS,Goerss08_QcohMfg,Lurie-chromatic,Peterson19_formalgeometry,Piotr-chromatic,Rognes-chromatic, Guillou-chromatic}.

In this section, we look at how $v_n$-periodic phenomena manifests in $\pi_{*,*}^F(\mathbb{S})$.

\subsection*{Conventions and Notation}
For (motivic or topological) spectra $\mathrm{E}$ and $\mathrm{X}$, we generally let $\mathrm{X}_\mathrm{E} = L_{\mathrm{E}}\mathrm{X}$ denote the Bousfield localization of $\mathrm{X}$ with respect to $\mathrm{E}$. There are a few exceptions to this convention which are clear from context. When $\mathrm{E}= \mathrm{H}\qq$ is the Eilenberg--Mac Lane spectrum associated to $\qq$, we write $\mathrm{X}_{\qq} = \mathrm{X}_{\mathrm{H}\qq}$ and call this the rationalization of $\mathrm{X}$. For $\mathrm{E} = \mathrm{H}\mathbb{F}_p$, we write $\mathrm{X}_p^\wedge = \mathrm{X}_{\mathrm{H}\mathbb{F}_p}$ and call this the $p$-completion of $\mathrm{X}$.

\subsection{Rationalization}
By Serre's finiteness theorem, the rationalization of the sphere $\ss^{\mathrm{top}}_\qq$ has very simple homotopy groups:
\[
\pi_s(\mathbb{S}^{\mathrm{top}}_\mathbb{Q}) \cong \left\{ \begin{array}{rl}
    \mathbb{Q} & s=0 \\
    0 & \mathrm{else.}
\end{array} \right.
\]
This simplicity extends to the entire stable homotopy category. There is an equivalence of spectra $\ss^{\mathrm{top}}_\mathbb{Q} \simeq \mathrm{H}\qq^{\mathrm{top}}$, so that
\[
\mathrm{Sp}_\qq \simeq \mathrm{Mod}_{\mathrm{H}\qq^{\mathrm{top}}}(\mathrm{Sp}) \simeq\euscr{D}(\qq),
\]
where $\euscr{D}(\qq)$ is the derived category of rational chain complexes. In this way, one often says that rational stable homotopy theory is equivalent to rational linear algebra.

Rationalization in motivic homotopy theory is more complicated. It was first observed by Morel that inverting $2 \in \pi_{0,0}^F(\mathbb{S})$ induces a splitting
\[
\mathrm{SH}(F)[1/2] \simeq \mathrm{SH}(F)^+ \times \mathrm{SH}(F)^-
\]
into a ``positive" and ``negative" part. As a result, this leads to a splitting of the rationalization
\[
\mathrm{SH}(F)_{\qq} \simeq \mathrm{SH}(F)_{\qq}^+ \times \mathrm{SH}(F)_{\qq}^-.
\]

The positive part $\mathrm{SH}(F)_{\qq}^+$ is more familiar. By work of Cisinski and Deglise, there is an equivalence of motivic spectra $\ss_{\qq}^+ \simeq \mathrm{H}\qq$, where $\mathrm{H}\qq$ is the motivic spectrum representing motivic cohomology with rational coefficients \cite{CisDeg_triangulatemixedmotives}. This leads to an identification of the positive part of the rational category
\[
\mathrm{SH}(F)_{\qq}^+ \simeq \mathrm{Mod}_{\mathrm{H}\qq}(\mathrm{SH}(F)) \simeq \euscr{DM}(\qq),\footnote{More generally, by work of R\"ondigs and \O stv\ae r \cite{RonOst06_motives}, and later of Elmanto and Kolderup \cite{ElmKol20_modulesovermotivicringspectra}, there is an equivalence $\mathrm{Mod}_{\mathrm{H}\mathbb{Z}}(\mathrm{SH}(F)) \simeq \euscr{DM}$, Voevodsky's category of motives \cite{Voe00_triangulatedcategoriesofmotives}.}
\]
where $\euscr{DM}(\qq)$ is Voevodsky's category of rational motives over $F$. 

The negative part $\mathrm{SH}(F)_{\qq}^-$ is less familiar. By work of Ananyevskiy, Levine, and Panin, the negative part vanishes precisely when $-1$ is a sum of squares in $F$. Equivalently, this occurs whenever $\mathrm{W}(F)$ contains a torsion free summand. For instance, $\mathrm{SH}(\mathbb{C})_{\mathbb{Q}}^{-}$ is trivial while $\mathrm{SH}(\mathbb{R})^-_{\mathbb{Q}}$ is not.
The authors show that there is an identification akin to the positive case
\[
\mathrm{SH}(F)_{\qq}^{-} \simeq \euscr{DM}_\mathrm{W}(\mathbb{Q}),\footnote{One can also $\mathrm{SH}(F)_{\mathbb{Q}}^-$ as $\mathrm{Mod}_{\mathrm{HW}_{\mathbb{Q}}}(\mathrm{SH}(F))$, where $\mathrm{HW}$ is the Witt cohomology spectrum. This spectrum may be identified as the cofiber of the periodicity operator on $\mathrm{kw}$.}
\]
where $\euscr{DM}_\mathrm{W}(\mathbb{Q})$ is a category of rational Witt motives over $\mathrm{F}$. Using this identification, they show an analogue of Serre's Finiteness Theorem for $\mathbb{S}_{\mathbb{Q}}^-$.
\begin{thm}[\cite{AnaLevPan17_Wittsheaves_etainverted}]
    For $n>0$, the group $\pi_{n+k, k}^F(\ss_{\qq}^-)$ vanishes.
\end{thm}

We will not discuss rational structures further here. See the work of Garkusha and of D\'eglise, Fasel, Jin, and Khan for further studies \cite{Gar19_rationalmotivic, DegFasJinKha21_rationalmotivic}.

\subsection{$v_1$-periodic phenomena}
\label{subsection: v1 periodicity}
Away from rationalization, chromatic homotopy theory is best understood at height $n=1$. 
There are two computational tools we will focus on for computing $v_1$-periodic homotopy. For the rest of this section, all spectra are implicitly $2$-complete, and all motivic spectra are implicitly $(2,\eta)$-complete.

\subsubsection{Topological $v_1$-periodicity}
There is a cofiber sequence
\[
\ss^{\mathrm{top}}_{\mathrm{K}(1)} \to \mathrm{KO} \xrightarrow{\psi^3-1} \mathrm{KO},
\]
where $\mathrm{K}(1) = \mathrm{KU}/2$ is the height 1 Morava K-theory spectrum. This expresses the $\mathrm{K}(1)$-local sphere as the equalizer of the Adams operations and the identity map on $\mathrm{KO}$. The homotopy groups $\pi_*(\mathrm{KO})$ are well known:
\[
\pi_*(\mathrm{KO}) \cong \frac{\mathbb{Z}[\eta, \alpha, \beta^{\pm 1}]}{(2\eta, \eta^3, \eta\alpha, \alpha^2-4\beta)}, \quad |\eta| =1,\, |\alpha|=4, \, |\beta| = 8.
\]

The homotopy groups $\pi_{*}(\ss^{\mathrm{top}}_{\mathrm{K}(1)})$ are closely related to the $v_1$-periodic elements in $\pi_{*}(\mathbb{S}^{\mathrm{top}})$. However, it accounts for much more. For example, there is a nonzero element $\zeta \in \pi_{-1}(\mathbb{S}^{\mathrm{top}}_{\mathrm{K}(1)})$ coming from the map $\partial$ in the long exact sequence:
\[
\cdots \to\mathbb{Z} \cong\pi_0(\mathrm{KO}) \xrightarrow{\partial} \pi_{-1}(\mathbb{S}^{\mathrm{top}}_{\mathrm{K}(1)}) \to \pi_{-1}(\mathrm{KO}) \cong 0 \to \dots
\]
Since $\ss^{\mathrm{top}}$ is connective, this class $\zeta$ cannot represent a $v_1$-periodic element in its stable homotopy. Even worse, one can show that $\ss_{\mathrm{K}(1)}^\mathrm{top}$ is not even bounded below.

One way to bypass this complication is to instead work with a cofiber sequence where all terms are bounded below. Let $\mathrm{ko}$ denote the connective cover of $\mathrm{KO}$ and let $\mathrm{ksp}$ denote its 4-connective cover, so that $\mathrm{ksp}[4] \simeq \tau_{\geq 4}\mathrm{KO}.$ The map $\psi^3-1$ lift to a map on $\mathrm{ko}$ which factors through $\mathrm{ksp}$, giving a cofiber sequence
\[
\mathrm{j} \to \mathrm{ko} \xrightarrow{\psi^3-1} \mathrm{ksp}[4]
\]
whose fiber we denote by $\mathrm{j}$. The homotopy groups $\pi_*(\mathrm{j})$ are easily determined by a long exact sequence in homotopy groups.

Adams showed that the Hurewicz map $\mathbb{S} \to \mathrm{j}$ is surjective on homotopy groups \cite{Adams-JX-IV}. A key insight of Mahowald was that this map induces an isomorphism on $v_1$-periodic elements, hence determining the $v_1$-periodic elements in $\pi_*(\mathbb{S}^{\mathrm{top}})$. To show this, Mahowald made extensive computations with the $\mathrm{ko}$-based Adams spectral sequence, which is referred to as the $\mathrm{bo}$-resolution \cite{Mah81}.\footnote{Although we use $\mathrm{ko}$ to refer to connective real K-theory, we will still refer to this spectral sequence as the $\mathrm{bo}$-resolution.} 

\subsubsection{Motivic $v_1$-periodicity}
Upon inspection of motivic Adams charts one begins to notice lifts of topological $v_1$-periodic elements to $\pi^F_{*,*}(\mathbb{S})$. Recall the hermitian K-theory spectrum $\mathrm{KQ}$, its very effective cover $\mathrm{kq}$, and the algebraic K-theory spectrum $\mathrm{KGL}$. These spectra often behave as motivic analogues of $\mathrm{KO}$, $\mathrm{ko}$, and $\mathrm{KU}$, respectively. Taking inspiration from the topological cofiber sequence, Balderamma, Ormsby, and Quigley have shown the following.
\begin{thm}[\cite{BOQ23}]
\label{thm:BOQ}
    Let $F$ be any field such that $\mathrm{char}(F) \neq 2$, and let $\mathrm{X} \in \mathrm{SH}(F)$ be any motivic spectrum. There is a cofiber sequence
    \[
    \mathrm{X}_{\mathrm{KGL}/2} \to (\mathrm{KQ} \otimes \mathrm{X})_{2, \eta}^\wedge \xrightarrow{\psi^3-1} (\mathrm{KQ} \otimes \mathrm{X})_{2, \eta}^\wedge.\footnote{We emphasize the need for completion after taking the tensor product.}
    \]
\end{thm}

\begin{example} 
    The hermitian K-theory of $\mathbb{C}$ is well known, and leads to an identification of the homotopy of $\mathrm{KQ}$:
    \[
    \pi_{*,*}^{\mathbb{C}}(\mathrm{KQ}) \cong \frac{\mathbb{Z}[\eta, \alpha, v_1^{\pm 4}, \tau]}{(2\eta, \tau\eta^3, \eta \alpha, \alpha^2-4v_1^4)}.
    \]
    Upon Betti realization, the motivic Adams operation $\psi^3-1:\mathrm{KQ} \to \mathrm{KQ}$ is sent to the topological Adams operation $\psi^3-1:\mathrm{KO} \to \mathrm{KO}$. Thus, 
    \[
    (\psi^{3}-1)(\eta) = 0, \quad (\psi^3-1)(\alpha) = 3^2\alpha, \quad (\psi^3-1)(v_1^{4}) = 3^4v_1^4, \quad (\psi^3-1)(\tau) = 0.
    \]
    We display the result in \Cref{fig:kgl/2 sphere over c} in $(s,w)$-charts. With our convention, after determining the action $\psi^3-1$, this chart also depicts $\pi_{*,*}^{\mathbb{C}}(\ss_{\mathrm{KGL}/2})$.
    \begin{figure}
        \centering
        \includegraphics[scale=.75]{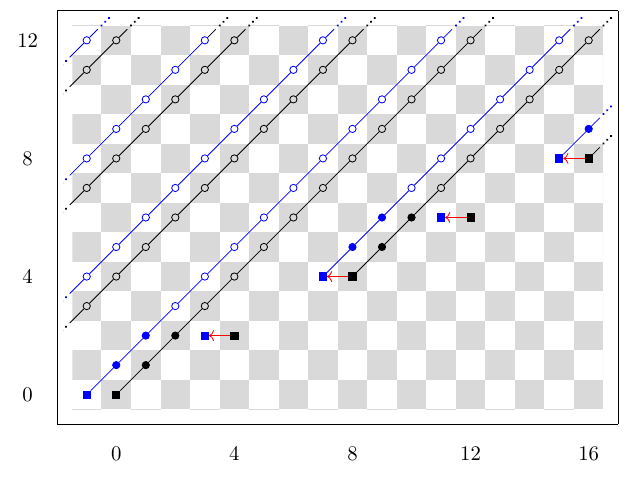}
        \caption{The map $\psi^3-1:\mathrm{KQ} \to \mathrm{KQ}$. A $\blacksquare$ indicated $\mathbb{Z}[\tau]$, a $\bullet$ indicates $\mathbb{F}_2[\tau]$, and a $\circ$ indicates $\mathbb{F}_2$. The blue and black colors indicate which copy of $\mathrm{KQ}$ a class belongs to. Red arrows indicate nonzero values for $\psi^3-1$.}
        \label{fig:kgl/2 sphere over c}
    \end{figure}
\end{example}

By Morel's vanishing result, $\pi_{n+k,n}^F(\mathbb{S}) = 0$ for $k < 0$. However, the topological element $\zeta \in \pi_{-1}(\ss^{\mathrm{top}}_{\mathrm{K}(1)})$ lifts through the cofiber sequence of \Cref{thm:BOQ} to produce a nonzero element $\zeta \in \pi_{-1, 0}^F(\mathbb{S}_{\mathrm{KGL}/2})$. Thus the homotopy of $\mathbb{S}_{\mathrm{KGL}/2}$ will account for more than just the $v_1$-periodicity in $\mathbb{S}$.

As in the topological case, one may remedy this by passing to covers. Let $\mathrm{ksp}$ denote the very effective cover of $\mathrm{KQ}[4,2]$. By \cite[Corollary 3.32]{BH21}, the Adams operations on $\mathrm{KQ}$ lift to a map on the very effective cover $\mathrm{kq}$ which factors through $\mathrm{ksp}$, giving a cofiber sequence
\[
\mathrm{jq} \to \mathrm{kq} \xrightarrow{\psi^3-1} \mathrm{ksp}[4,2].
\]
The homotopy of $\mathrm{jq}$ is easily determined by a long exact sequence in homotopy groups. One way to depict its homotopy groups is to take those given by \Cref{fig:kgl/2 sphere over c} and delete all classes which lie above the line $w = s$.

Inspired by Mahowald's findings, one expects that the $v_1$-periodic elements of $\pi_{*,*}^F(\mathbb{S})$ are entirely captured by $\pi_{*,*}^F(\mathrm{jq})$. One might further hope to show this by computing the $\mathrm{kq}$-based Adams spectral sequence, known as the $\mathrm{kq}$-\emph{resolution}. This spectral sequence takes the form
\[
\mathrm{E}_1^{s,f,w} = \pi_{s+f, w}^F(\mathrm{kq} \otimes \overline{\mathrm{kq}}^{\otimes f}) \implies \pi_{s,w}^F(\mathrm{kq}).
\]
As $\mathrm{kq}$ is not flat over $\mathrm{kq} \otimes \mathrm{kq}$, the $\mathrm{E}_1$-page of the $\mathrm{kq}$-resolution does not identify with the cobar complex computing the cohomology of the Hopf algebroid $\pi_{*,*}^F(\mathrm{kq} \otimes \mathrm{kq})$. Therefore, understanding the $\mathrm{E}_1$-page of the $\mathrm{kq}$-resolution is of key importance, and indeed requires substantial computation.

The study of the $\mathrm{kq}$-resolution was initiated by Culver and Quigley. We summarize their findings in the following.
\begin{thm}[\cite{CQ21}]
The $v_1$-periodic elements of $\pi_{*,*}^\mathbb{C}(\mathbb{S})$ are concentrated in the $0$ and $1$-lines of the $\mathrm{kq}$-resolution. Furthermore, there is a summand $\pi_{*,*}^\mathbb{C}(\mathrm{ksp}) \subseteq \pi_{*,*}^\mathbb{C}(\mathrm{kq} \otimes \overline{\mathrm{kq}})$, and the $d_1$-differential
\[
d_1:\pi_{*,*}^\mathbb{C}(\mathrm{kq}) \to \pi_{*,*}^\mathbb{C}(\mathrm{ksp}[4,2]) \subseteq \pi_{*,*}^\mathbb{C}(\mathrm{kq}\otimes \overline{\mathrm{kq}})
\]
restricts to this summand and is given by $\psi^3-1$. There are no other differentials from the 0-line and all elements of the complementary summand of the 1-line all support differentials. Thus, the $v_1$-torsion free component of $\pi_{*,*}^\mathbb{C}(\mathrm{kq})$ are isomorphic to $\pi_{*,*}^\mathbb{C}(\mathrm{jq})$. 
\end{thm}

Recent progress has been made by the author, Petersen, and Tatum to understand the $\mathrm{E}_1$-page of the $\mathrm{kq}$-resolution over the base fields $F\in\{\mathbb{R}, \mathbb{F}_p, \mathbb{Q}_p, \mathbb{Q}\}$ \cite{Realkqcoop,finitekqcoop,MorPetTat}.

Closely related to $\mathrm{jq}$ is the spectrum $\mathrm{L} = \mathrm{fib}(\mathrm{kq} \xrightarrow{\psi^3-1}\mathrm{kq})$. From many points of view, this spectrum captures nearly as much $v_1$-periodicity as $\mathrm{jq}$ and is much more computationally accessible. The homotopy of $\mathrm{L}$ has been computed via the effective slice spectral sequence for $F=\mathbb{C}$ and $\mathbb{R}$ by Belmont, Isaksen, and Kong \cite{belmontisaksenkong-v1R}, and for $F=\mathbb{F}_q, \mathbb{Q}_p,$ and $\mathbb{Q}$ by Kong and Quigley \cite{kongquigley}. 

There has also been extensive study of $v_1$-periodic elements over a variety of fields by Quigley's work on the motivic Mahowald invariant; see \cite{Quigley19, Qui21_RC2mahowald, Qui21_generalmahowald}.

\begin{remark}
    Mahowald's work on the $\mathrm{bo}$-resolution may be interpreted as a solution to the height 1 telescope conjecture at the prime 2. One may then hope to attack a motivic analogue by using the $\mathrm{kq}$-resolution. However, exotic periodicity in motivic homotopy theory makes even the statement of a motivic telescope conjecture more opaque (see also \cite[Section 6]{CQ21}). Apart from an increase in difficulty in computations, there are significant structural differences. For example, there is a topological equivalence $L_{\mathrm{KO}}(-) \simeq L_{\mathrm{KU}}(-)$, which one can use to show that the $v_1$-localized $\mathrm{bo}$-resolution converges to the $\mathrm{K}(1)$-localization. The analogous statement is not true in motivic homotopy. Let $\mathrm{KW}=\eta^{-1}\mathrm{KQ}$ be the \emph{Witt} $\mathrm{K}$-\emph{theory} spectrum. By \cite[Proposition 6.39]{CQ21} we have that $L_{\mathrm{KGL}}(-) \simeq L_{\mathrm{KQ}\oplus \mathrm{KW}}(-)$, making an algebraic comparison with the $v_1$-periodic $\mathrm{kq}$-resolution less straightforward.
\end{remark}

\subsection{$v_2$-periodic phenomena}
Much of the computational work in chromatic homotopy theory over the last 30 years has been conducted at height 2. Analogous to the role that $\mathrm{ko}$ plays in detecting $v_1$-periodicity, the connective spectrum $\mathrm{tmf}$ of \emph{topological modular forms} plays an important role in detecting $v_2$-periodicity \cite{BehHilHopMah09_v2_32, BehManQui23_tmfhurewiczimage, cardav25_periodicfamilies,BobQui26_etafamilies,BBQ26}. In motivic homotopy theory, the study of $v_2$-periodicity is in its nascent stages.

Over $F=\mathbb{C}$, Gheorghe, Isaksen, Krause, and Ricka used the filtered analogue functor $\Gamma_\star:\mathrm{Sp} \to \mathrm{Fil}(\mathrm{Sp})$ to construct a \emph{motivic modular forms spectrum} $\mathrm{mmf} = \Gamma_\star(\mathrm{tmf}) \in \mathrm{SH}(\mathbb{C})$ \cite{GIKR22}. The homotopy groups $\pi_{*,*}^\mathbb{C}(\mathrm{mmf})$ have been computed by Isaksen, Kong, Li, Ruan, and Zhu by both the $\textbf{mASS}^\mathbb{C}(\mathrm{mmf})$ \cite{isaksen_ASS_mmf} and the $\textbf{mANSS}^\mathbb{C}(\mathrm{mmf})$ \cite{IKLRZ24_ANSSmmf}. 

The usage of $\mathrm{mmf}$ in computations in $\mathbb{C}$-motivic homotopy theory is multifaceted. The unit map $\mathbb{S} \to \mathrm{mmf}$ induces a map of spectral sequences $\textbf{mASS}^\mathbb{C}(\mathbb{S}) \to \textbf{mASS}^\mathbb{C}(\mathrm{mmf})$. This allows one to pull back differentials for $\mathrm{mmf}$, whose Adams spectral sequence is completely understood, to differentials for $\mathbb{S}$. This is a powerful technique leveraged by Isaksen, Wang, and Xu \cite{IWX23}. Additionally, the elements in the Hurewicz image of $\pi_{*,*}^\mathbb{C}(\mathbb{S}) \to \pi_{*,*}^\mathbb{C}(\mathrm{mmf})$ are necessarily $\Delta^8 = v_2^{32}$-periodic. The previous authors have computed the image of this map up to stem 90, identifying many $v_2$-periodic families in $\pi_{*,*}^\mathbb{C}(\mathbb{S})$.

The construction of a motivic modular forms spectrum has only been performed over $\mathrm{F}=\mathbb{C}$,\footnote{Relatedly, Carrick, Davies, and van Nigtvecht have studied $\mathrm{BP}$-synthetic analogues (i.e. $\mathbb{C}$-motivic) of $\mathrm{TMF}$ and $\mathrm{Tmf}$ \cite{cardav25_periodicfamilies,CarDavNig25_descentssviasynthetic}.} and so the above techniques have as of yet not been extended to other base fields. However, assuming the existence of such a spectrum $\mathrm{mmf}$ in $\mathrm{SH}(\mathbb{R})$, the $\mathrm{E}_2$-page of the $\textbf{mASS}^\mathbb{R}(\mathrm{mmf})$ must take the form\footnote{As a comodule algebra over the dual Stenrod algebra, $\mathrm{H}_*(\mathrm{tmf}) \cong (\aeu // \aeu (2))^\vee$. For any motivic analogue, we expect the analogous isomorphism to hold over the dual Steenrod algebra.}
\[
\mathrm{E}_2^{s,f,w} = \mathrm{Ext}^{s,f,w}_{\aeu(2)^\vee}(\mathbb{M}_2^\mathbb{R}, \mathbb{M}_2^\mathbb{R}) \implies \pi_{s,w}^\mathbb{R}(\mathrm{mmf}).
\]
Recent work of Emming has computed this $\mathrm{E}_2$-page, and assuming the existence of a motivic modular forms spectrum, computes many of the $d_2$-differentials in the $\textbf{mASS}^\mathbb{R}(\mathrm{mmf})$ \cite{emming26_ExtA2overR}.

\begin{remark}
    The concept of evenness has been an ever present piece of chromatic homotopy theory. For example, $\pi_*(\mathrm{MU})$ is concentrated in even degrees, as is any complex oriented spectrum. It has become more present over the last 10 years that this evenness is in fact very powerful. For example, the even filtration of Hahn, Raksit, and Wilson has given new ways to study invariants of ring spectra by first resolving them by even covers \cite{HRW}. This filtration has been shown to agree with the decalage of the Adams--Novikov spectral sequence. By the work of Gheorghe, Isaksen, Krause, and Ricka, this shows that one can model $\mathrm{SH}(\mathbb{C})^{\mathrm{cell}}$ simply by looking at filtered spectra which are modules over $\mathbb{S}^{\mathrm{ev}}$, the even filtration applied to the sphere.

    The filtered analogue functor also behaves well with respect to evenness in another way. Namely, if $\mathrm{X} \in \mathrm{Sp}$ has even $\mathrm{BP}$-homology, then the filtered analogue $\Gamma_\star(\mathrm{X})$ is particularly simple to work with. This has been particularly useful for studying periodicity in motivic homotopy theory thus far. Forthcoming work of Emming also shows that, for $\mathrm{X}$ connective, having even $\mathrm{BP}$-homology ensures a tight connection between the $\textbf{ANSS}(\mathrm{X})$ and the $\textbf{SliceSS}^\mathbb{C}(\Gamma_\star(\mathrm{X}))$

    To be precise, there are a class of spectra known as \emph{connective higher real $\mathrm{K}$-theories} which are very useful for studying periodicity. These arise as particular connective variants of the spectra $\mathrm{eo}_n(\mathrm{G}) = \mathrm{E}_n^{hG}$, where $\mathrm{E}_n$ is the height $n$ Lubin--Tate spectrum and $G$ is a finite subgroup of the extended Morava Stabilizer group $\mathbb{G}_n$. The connective K-theory $\mathrm{ko}$ is a model for $\mathrm{eo}_1(C_2)$, and various connective topological modular forms spectra (specifically $\mathrm{tmf}_0(3)$ and $\mathrm{tmf}$) are related to various $\mathrm{eo}_2(G)$ (specifically $\mathrm{eo}_2(C_2)$ and $\mathrm{eo}_2(G_{24})$). These spectra all have even $\mathrm{BP}$-homology, and their analogues $\mathrm{kq} = \Gamma_\star(\mathrm{ko})$ and $\mathrm{mmf} = \Gamma_\star(\mathrm{tmf})$ have proven themselves very useful in motivic homotopy theory.

    Unfortunately, our luck runs out. Recent work of Carrick and Hill shows that for $|G|>2$, the $\mathrm{BP}$-homology of $\mathrm{eo}_n(G)$ is nonzero in odd degrees, and for $G=C_2$ the same as true for $n \geq 3$ \cite{CarHil25_MUhomology_eonG}. In particular, the filtered analogue of a higher real K-theory is quite difficult to work with. Thus, alternative methods must be invoked to study periodicity in motivic homotopy theory.
\end{remark}

\section{Exotic periodic phenomena}
\label{section:exotic}
There is interesting periodic phenomena in motivic homotopy which is absent in topological stable homotopy theory. For example, we have seen that $\eta \in \pi_{1,1}^F(\mathbb{S})$ is non-nilpotent. This leads one to consider $\eta$-periodic phenomena which is unique to motivic homotopy theory.\footnote{Unique in that it does not appear topologically. In $C_2$-equivariant stable homtoopy theory, the Hopf map $\eta \in \pi_\sigma^{C_2}(\mathbb{S})$ is also non-nilpotent.}

In this section, we investigate exotic periodic phenomena in motivic homotopy theory stemming from this observation.

\subsection{$\eta$-periodic phenomena}

As the first Hopf map $\eta \in \pi_{1,1}^F(\ss)$ is non-nilpotent, one may form the telescope
\[
\eta^{-1}\mathbb{S} = \mathrm{colim}(\ss \xrightarrow{\eta}\ss[-1,-1] \xrightarrow{\eta}\ss[-2,-2]\xrightarrow{\eta}\cdots),
\]
which we will call the $\eta$-\emph{periodic sphere}. In parallel to the simplicity of $\pi_*(\ss^{\mathrm{top}}_\mathbb{Q})$, one would expect that the $\eta$-periodic sphere has relatively simple homotopy. 

The computation of the homotopy of $\eta^{-1}\ss$ was first approached by Guillou and Isaksen over $F=\mathbb{C}$ \cite{GuiIsa15}. Their method of approach was to use the motivic Adams spectral sequence. As $\eta$ is detected by the element $h_1 \in \mathrm{Ext}_{\euscr{A}^\vee}^{\cc}(\ss)$, they deduce that $h_1$ must be non-nilpotent, and hence the $\textbf{mASS}^\mathbb{C}(\eta^{-1}\mathbb{S})$ takes the form
\[
\mathrm{E}_2^{s,f,w} = h_1^{-1}\mathrm{Ext}_{\euscr{A}^\vee}^{s,f,w}(\mathbb{C}) \implies \pi_{s,w}^\mathbb{C}(\eta^{-1}\mathbb{S}).
\]
Their main result is a computation of this $\mathrm{E}_2$-page.
\begin{thm}[{\cite[Theorem 1.1]{GuiIsa15}}]
    The $\mathrm{E}_2$-page of the $\textup{\textbf{mASS}}^{\mathbb{C}}(\eta^{-1}\mathbb{S})$ takes the form
    \[
    h_1^{-1}\mathrm{Ext}_{\euscr{A}^\vee}^{\cc}(\ss) \cong \mathbb{F}_2[h_1^{\pm 1}, P,v_2, v_3, \dots],
    \]
    where $|h_1| = (1,1,1)$, $|P| = (8,4,4)$, and $|v_n| = (2(2^n-1), 1, 2^n-1).$
\end{thm}

While Guillou and Isaksen do not compute the homotopy of $\eta^{-1}\ss$, they conjecture that there are differentials in the $\textbf{mASS}^{\mathbb{C}}(\eta^{-1}\ss)$ for all $k \geq 0$ and $n \geq 3$:
\[
d_2(P^kv_n) = P^kv_{n-1}^2.
\]
For degree reasons, this would imply that the $\textbf{mASS}^\mathbb{C}(\eta^{-1}\ss)$ collapses on the $\mathrm{E}_3$-page.

\begin{remark}
    The $\mathrm{E}_2$-page of the $\textbf{mASS}^\mathbb{C}(\eta^{-1}\mathbb{S})$ offers us a glimpse into the peculiar relationship that classical and exotic periodicity have with each other.

    In topology, of the periodicity classes $v_n$, the $\mathrm{E}_2$-page of the $\textbf{ASS}(\ss^{\mathrm{top}})$, i.e. the cohomology $\mathrm{Ext}_{\aeu^\vee_{\mathrm{top}}}(\ss^\mathrm{top})$, contains only $v_0 = 2$ and no other $v_n$'s. Inverting $v_0$ on the $\mathrm{E}_2$-page yields $\mathbb{F}_2[v_0^{\pm 1}]$. To see higher $v_n$'s, we instead pass to the quotient.
    
    Consider the mod-2 Moore spectrum $\ss^{\mathrm{top}}/2 = \mathrm{cofib}(\mathbb{S}^{\mathrm{top}} \xrightarrow{2} \mathbb{S}^{\mathrm{top}})$. The Periodicity Theorem of Devinatz, Hopkins, and Smith ensures that there is some $v_1$-self-map on this cofiber \cite{DHS88,HS98}. Adams showed that a minimal self-map is of periodicity 4, giving a non-nilpotent self-map $v_1^4:\mathbb{S}^{\mathrm{top}}/2[8] \to \mathbb{S}^{\mathrm{top}}/2$ \cite{Adams-JX-IV}.  

    The composition $\ss^{\mathrm{top}}[8] \to \ss^{\mathrm{top}}/2[8] \xrightarrow{v_1^4} \ss^{\mathrm{top}}/2$
    is essential, where the first map includes the bottom cell of the Moore spectrum. It follows that there is some non-nilpotent element of $\pi_8(\ss^{\mathrm{top}}/2)$ witnessing $v_1^4$. Using an algebraic Atiyah--Hirzebruch spectral sequence
    \[
    \mathrm{Ext}_{\euscr{A}^\vee_\mathrm{top}}(\ss^{\mathrm{top}})\{x_0, x_1\} \implies \mathrm{Ext}_{\euscr{A}^\vee_\mathrm{top}}(\ss^{\mathrm{top}}/2)
    \]
    one discovers that there is a unique class $v_1^4=P$ of degree $(8,4)$ on the $\mathrm{E}_2$-page of the $\textbf{ASS}(\ss^{\mathrm{top}}/2)$ which must witness this self-map. Thus $P$ must be non-nilpotent and survive the $\textbf{ASS}(\mathbb{S}^{\mathrm{top}}/2)$. More generally, if $\mathrm{X}$ is a type $n$ spectrum, then there will be some power of $v_n$ witnessed in the $\mathrm{E}_2$-page of the $\textbf{ASS}(\mathrm{X})$ which is non-nilpotent.

    The computation of Guillou and Isaksen implies that inverting an exotic periodic element is significantly more interesting. Indeed, as soon as one inverts $h_1$, each $v_n$ becomes visible and non-nilpotent at the level of motivic Adams $\mathrm{E}_2$-pages.
\end{remark}

Andrews and Miller approached the $\eta$-periodic sphere over $F=\mathbb{C}$ via the motivic Adams--Novikov spectral sequence \cite{AndMil17}. As $\eta$ is detected by  $\alpha_1 \in \mathrm{Ext}^{1,1,1}_{\mathrm{BPGL}}(\mathbb{S})$, they deduce that $\alpha_1$ must be non-nilpotent, and hence the $\textbf{mANSS}^\mathbb{C}(\eta^{-1}\mathbb{S})$ takes the form
\[
\mathrm{E}_2^{s,f,w} = \alpha_1^{-1}\mathrm{Ext}^{\cc}_{\mathrm{BPGL}}(\mathbb{S}) \implies \pi_{s,w}^\mathbb{C}(\eta^{-1}\mathbb{S}).
\]
Unlike the topological $h_1$, the topological $\alpha_1$ is non-nilpotent on the $\mathrm{E}_2$-page of the $\textbf{ANSS}(\mathbb{S}^{\mathrm{top}})$. The key result of Andrews--Miller is a computation of this localized Ext group. Using the algebraic Novikov spectral sequence, they show that there is an isomorphism
\[
\alpha_1^{-1}\mathrm{Ext}^{*,*}_{\mathrm{BP}_*\mathrm{BP}} \cong \mathbb{F}_2[\alpha_1^{\pm 1}, \alpha_3, \alpha_4]/(\alpha_4^2).
\]
As $\eta^{-1}\mathbb{S}^{\mathrm{top}} \simeq *$, the $\textbf{ANSS}(\eta^{-1}\mathbb{S}^{\mathrm{top}})$ sequence must converge to zero. The topological differential $d_3(\alpha_4) = \alpha_1^4$ persists to the localization, and this determines the spectral sequence. 

Betti realization lifts the topological result to the $\mathbb{C}$-motivic setting.
\begin{thm}[\cite{AndMil17}]
    The $\mathrm{E}_2$-page of the $\textup{\textbf{mANSS}}^{\mathbb{C}}(\eta^{-1}\mathbb{S})$ takes the form
    \[
    \alpha_1^{-1}\mathrm{Ext}^{\cc}_{\mathrm{BPGL}}(\mathbb{S}) \cong \mathbb{F}_2[\alpha_1^{\pm 1}, \alpha_3, \alpha_4, \tau]/(\alpha_4^2).
    \]
    The spectral sequence is determined via the Liebniz rule by the differential $d_3(\alpha_4) = \tau \alpha_1^4$. There are no hidden extensions on the $\mathrm{E}_\infty$-page, and so there is an isomorphism
    \[
    \pi_{*,*}^\mathbb{C}(\eta^{-1}\mathbb{S}) \cong \mathbb{F}_2[\eta^{\pm 1}, \mu_9, \epsilon]/(\epsilon^2). \footnote{As $\eta$ is $v_1$-periodic, all of the classes of positive stem degree of $\pi_{*,*}^{\mathbb{C}}(\eta^{-1}\mathbb{S})$ are detected in $\pi_{*,*}^{\mathbb{C}}(\mathrm{jq})$. Upon inspection, these are precisely the 2-torsion classes detected along the vanishing line of slope $\frac{1}{2}$ in the $\textbf{mASS}^\mathbb{C}(\mathbb{S})$. This observation is made more precise by Culver and Quigley in their work on the $\mathrm{kq}$-resolution \cite{CQ21}. Note that in terms of traditional Adams spectral sequence naming conventions, $\mu_9$ is detected by $Ph_1$ and $\epsilon$ is detected by $c_0$.
    }
    \]
\end{thm}
As a consequence, the conjectured differentials in the $\textbf{mASS}^{\mathbb{C}}(\eta^{-1}\mathbb{S})$ are indeed the only ones which occur. 

Using this computation, Guillou and Isaksen computed the $\mathbb{R}$-motivic analogue of their previous work \cite{GuiIsa16_Rmotivic_etaperiodic}. The $\rho$-Bockstein spectral sequence
\[
h_1^{-1}\mathrm{Ext}_{\euscr{A}^\vee}^{***}(\mathbb{C})[\rho] \cong \mathbb{F}_2[h_1^{\pm 1}, P, v_2, v_3, \dots, \rho] \implies h_1^{-1}\mathrm{Ext}_{\euscr{A}^\vee}^{***}(\mathbb{R})
\]
computes the $\mathrm{E}_2$-page of the $\textbf{mASS}^\mathbb{R}(\eta^{-1}\mathbb{S})$ from the $\mathbb{C}$-motivic analogue. They show that the only differentials are those generated under the Lieniz rule by $d_{2^{n}-1}(P^{n-1})=\rho^{2^{n}-1}v_n$. Base change from $\mathbb{R}$ to $\mathbb{C}$ determines all $d_2$-differentials, and an application of Moss's convergence theorem determines all higher differentials. 

Wilson extended the findings of Andrews and Miller to compute $\pi_{*,*}^{F}(\eta^{-1}\mathbb{S})$, where $F \in \{\mathbb{Q}, \mathbb{Q}_p\}$ \cite{Wil18_etainvertedrationals}. Wilson's computation proceeds in the same way as Guillou and Isaksen. The analogous $\rho$-Bockstein spectral sequence collapses over $F=\mathbb{Q}_p$ on a finite page as $\rho$ is nilpotent for these particular fields. Local computations over $F=\mathbb{Q}_p$ can then be used to deduce the result over $F=\mathbb{Q}$.

In the work of Guillou and Isaksen and of Wilson, interesting torsion appears in the homotopy groups of $\pi^F_{*,*}(\eta^{-1}\mathbb{S})$ which is absent in Andrews and Miller's computation over $\mathbb{C}$. Wilson observed that this torsion is related to the Witt ring $\mathrm{W}(F)$ of the base field. This observation was expanded on by Ormsby and R\"ondigs, who took a third approach to the $\eta$-periodic sphere \cite{OrmRon20}. Recall from \Cref{subsection: sliceSS} that the $\textbf{SliceSS}^F(\mathbb{S})$ takes the form
\[
\mathrm{E}_1^{s,q,w} = \pi_{s,w}^F(s_q(\mathbb{S})) \implies \pi_{s,w}^F(\mathbb{S}),
\]
and that the slices of the motivic sphere spectrum are determined by the topological Adams--Novikov $\mathrm{E}_2$-page. The class $\alpha_1 \in \mathrm{Ext}_{\mathrm{BP}_*\mathrm{BP}}^{*,*}$ which detects $\eta \in \pi_1(\mathbb{S}^\mathrm{top})$ acts by multiplication on the Adams--Novikov $\mathrm{E}_2$-page, hence induces a self-map of the slice spectral sequence $\mathrm{E}_1$-page.

\begin{thm}[\cite{OrmRon20}]
    Let $F$ be a field of finite cohomological dimension, $\mathrm{char}(F) \neq 2$. Then there is an equivalence of spectral sequences 
    \[
    \alpha_1^{-1}\textup{\textbf{SliceSS}}^F(\mathbb{S}) \simeq \textup{\textbf{SliceSS}}^F(\eta^{-1}\mathbb{S}),
    \]
    which take the form
    \[
    \mathrm{E}_1^{s,q,w} = \pi_{s,w}^F(\mathrm{H}\mathbb{F}_2[\alpha_1^{\pm 1}, \alpha_3, \alpha_4]/(\alpha_4^2)) \implies \pi_{s,w}^F(\eta^{-1}\mathbb{S}).
    \]
    Further, if $-1$ is a sum of 4 squares in $F$,\footnote{For example, this condition holds for $F \in \{\mathbb{C}, \mathbb{F}_p, \mathbb{Q}_p, \mathbb{Q}\}$ for $p$ odd, and does not hold for $F=\mathbb{R}$ or any formally real field \cite{Lam05_intro_quadraticforms}.} then
    \[
    \pi_{*,*}^F(\eta^{-1}\mathbb{S}) \cong \mathrm{W}(F)[\eta^{\pm 1}, \mu_9, \epsilon]/(\epsilon^2).
    \]
\end{thm}

Ormsby and R\"ondigs determine the differentials in the $\textbf{SliceSS}^F(\eta^{-1}\mathbb{S})$ by comparing with the $\textbf{SliceSS}^F(\mathrm{kw})$, where $\mathrm{kw} = \eta^{-1}\mathrm{kq}$ is the \emph{connective Witt} $\mathrm{K}$-\emph{theory}. The $\textbf{SliceSS}^F(\mathrm{kw})$ is completely understood by work of Ananyevskiy, R\"ondigs, and \O stv\ae r \cite{ARO20}, and one can pull back differentials by using that the unit map $\mathbb{S} \to \mathrm{kw}$ factors through $\eta^{-1}\mathbb{S}$. 

Wilson has indicated that the homotopy of $\eta^{-1}\mathbb{S}$ is seemingly related to Witt theory. As $\pi_{*,*}^F(\mathrm{kw}) \cong \mathrm{W}(F)[\beta^{\pm 1}]$, the use of $\mathrm{kw}$ by R\"ondigs and Ormsby is not surprising. However, the relationship is much deeper it appears at first glance.

The definitive study of $\eta$-periodic phenomena is due to Bachmann and Hopkins \cite{BH21}. While previous work had compared $\eta$-periodic motivic phenomena and $2$-periodic topological phenomena, Bachmann and Hopkins instead focus on the relationship with $v_1$-periodic topological phenomena.

\begin{thm}
    Let $F$ be a field, $\mathrm{char}(F) \neq 2$. Then there exists a cofiber sequence
    \[
    \eta^{-1}\mathbb{S} \to \mathrm{kw} \xrightarrow{\psi^3-1} \mathrm{kw}[4,2].
    \]
    As a consequence,
    \[
    \pi_{s,w}^F(\eta^{-1}\mathbb{S}) \cong \left\{ \begin{array}{ll}
        \mathrm{W}(F) & (s,w) = (n,n),\, n \in \mathbb{Z}\\
        \mathrm{coker}(8k:\mathrm{W}(F) \to \mathrm{W}(F)) & (s,w) = (n+k, n), \, k \equiv 3 \, (4), \, k > 0 \\
        \mathrm{ker}(8k: \mathrm{W}(F) \to \mathrm{W}(F)) & (s,w) = (n+k, n), \, k \equiv 0 \, (4), \, k >0 \\
        0 & \mathrm{else.}
    \end{array} \right.
    \]
\end{thm}

The map $\mathrm{kw} \xrightarrow{\psi^3-1} \mathrm{kw}[4,2]$ is, as the name suggests, derived from the Adams operations used to study $v_1$-periodicity in \Cref{subsection: v1 periodicity}, and the cofiber sequence for $\eta^{-1}\mathbb{S}$ essentially amounts to inverting $\eta$ on the cofiber sequence $\mathrm{jq} \to \mathrm{kq} \to \mathrm{ksp}[4,2]$.

\begin{remark}
    The homotopy groups $\pi_{*,*}^F(\eta^{-1}\ss)$ need not be purely torsion. By the formula of Bachmann and Hopkins, this occurs precisely when $\mathrm{W}(F)$ contains a torsion free summand. As a consequence, whenever $\mathrm{W}(F)$ contains a torsion free summand, the rationalization $\pi_{*,*}^F((\eta^{-1}\mathbb{S})_\mathbb{Q})$ does not vanish. Ananyevskiy, Levine, and Panin show that there is an equivalence $(\eta^{-1}\ss_{\mathbb{Q}}) \simeq \ss_{\mathbb{Q}}^{-}$ \cite{AnaLevPan17_Wittsheaves_etainverted}.\footnote{In related work, Bachmann has shown that $(\eta^{-1}\mathbb{S})_{\mathbb{Q}} \simeq (\rho^{-1}\mathbb{S})_{\mathbb{Q}}$ \cite{Bachmann-real}.} Thus another interpretation to the negative part of the rational motivic homotopy category is
    \[
    \mathrm{SH}(F)_{\qq}^{-} \simeq \mathrm{Mod}_{(\eta^{-1}\ss)_{\qq}}(\mathrm{SH}(F)).
    \]
    If $\pi_{*,*}^F(\eta^{-1}\mathbb{S})$ is purely torsion, such as in the case of $F=\mathbb{C}$, this implies that $\mathbb{S}_{\mathbb{Q}}^- \simeq *$ and hence $\mathrm{SH}(F)_{\mathbb{Q}}^-$ is trivial.
\end{remark}

\begin{remark}
    Since $\eta$ is $v_1$-periodic, every $\eta$-periodic in $\pi_{*,*}^F(\mathbb{S})$ class is $v_1$-periodic. This suggests that one option to studying the $v_1$-periodicity is to decompose classes into two varieties: those which are $\eta$-periodic, which are handled by $\pi_{*,*}^F(\eta^{-1}\mathbb{S})$; and, those which are not $\eta$-periodic. This differs from the topological situation in that, as $\eta$ is nilpotent, there is no such decomposition of topological $v_1$-periodic classes. 
    
    Since $\eta$ is detected by $\mathrm{ko}$ and not by $\mathrm{ku}$, the $\mathrm{bo}$-resolution is better equipped to deal with $v_1$-periodicity than the $\mathrm{bu}$-resolution. However, the above decomposition of $v_1$-periodicity in motivic homotopy theory highlights a more distinguished role played by the $\mathrm{kgl}$-resolution, which is more easily computable than the $\mathrm{kq}$-resolution, in studying the $v_1$-periodic classes which are not $\eta$-periodic. This is one of the motivating ideas behind recent work on computing the $\mathrm{kgl}$-resolution \cite{LPT-kuRsplitting,MorPetTat-BPGL1, MorPetTat}.
\end{remark}

\subsection{$w_1$-periodic phenomena}
In this section, we work exclusively over $F=\mathbb{C}$ unless otherwise stated.

Miller suggests that non-nilpotence of $\eta$ is evidence for a new family of periodicity operators $w_n$, where $w_0 = \eta$ as $v_0 = 2$. Andrews shows this to be true \cite{Andrews18}. We briefly recall some background for motivation.

Recall the topological $v_1$-self-map $v_1^4:\mathbb{S}^\text{top}/2[8] \to \mathbb{S}^\text{top}/2.$ One can use this self-map to find infinite families in $\pi_*(\mathbb{S}^{\mathrm{top}})$ in the following manner. Let $\alpha$ be one of the stable elements $\{\eta, \eta^2, \eta^3, \epsilon, \epsilon\eta\}$. Then $\alpha \in \pi_d(\mathbb{S}^{\mathrm{top}})$ is 2-torsion, so we may lift through the boundary map $\pi_{d+1}(\mathbb{S}^{\mathrm{top}}/2) \xrightarrow{\partial} \pi_d(\mathbb{S}^{\mathrm{top}})$ to a class $\widetilde{\alpha} \in \pi_{d+1}(\ss^{\mathrm{top}}/2)$.

By suspending $\widetilde{\alpha}$ by $8k$, we end up with a stable map whose target is $\ss^{\mathrm{top}}[8k]$. To this we may apply $v_1^{4k}$, which is non-nilpotent, followed by the pinch map $\ss^{\mathrm{top}}/2 \xrightarrow{p}\ss^{\mathrm{top}}[1]$ onto the top cell of the Moore spectrum. The whole composite, which we call $v_1^{4k}(\alpha)$, takes the form
\[
v_1^{4k}(\alpha): \ss^{\mathrm{top}}[d+1+8k] \xrightarrow{\widetilde\alpha} \ss^{\mathrm{top}}/2[8k] \xrightarrow{v_1^{4k}} \ss^{\mathrm{top}}/2 \xrightarrow{p}\ss^{\mathrm{top}}[1].
\footnote{The composition $\ss^{\mathrm{top}}[8] \xrightarrow{\iota} \ss^{\mathrm{top}}/2[8] \xrightarrow{v_1^4} \ss^{\mathrm{top}}/2 \xrightarrow{p} \ss^{\mathrm{top}}[1]$ represents $8\sigma \in \pi_7(\mathbb{S}^{\mathrm{top}})$.
}
\]
The key idea of Adams is that for all values of $k \geq 1$, the classes $v_1^{4k}(\alpha) \in \pi_{8k+d}(\mathbb{S}^{\mathrm{top}})$ are essential \cite{Adams-JX-IV}. Since by construction $v_1^{4k}(\alpha)$ factors through a $v_1$-self-map, we call such a family a $v_1$-\emph{periodic family}. These are 2-torsion $v_1$-periodic elements in $\pi_*(\mathbb{S}^{\mathrm{top}})$, and by work of Mahowald, they are in fact the only 2-torsion $v_1$-periodic elements \cite{Mah81}.

The general philosophy layed out by Adams can be summarized as follows:
\begin{enumerate}
    \item Take a non-nilpotent element $x \in \pi_s(\mathbb{S}^{\mathrm{top}})$.
    \item Construct a non-nilpotent self map $v: \ss^{\mathrm{top}}/x[n] \to \ss^{\mathrm{top}}/x.$
    \item Find some $x$-torsion element $\alpha \in \pi_d(\mathbb{S}^{\mathrm{top}})$ and lift it through the boundary map to a class $\widetilde{\alpha} \in \pi_{d+1}(\ss^{\mathrm{top}}/x)$.
    \item See if the composite 
    \[
    \ss^{\mathrm{top}}[d+1+nk] \xrightarrow{\widetilde{\alpha}} \ss^{\mathrm{top}}/x[nk] \xrightarrow{v^n} \ss^{\mathrm{top}}/x \xrightarrow{p} \ss^{\mathrm{top}}[s+1]
    \]
    is essential.
\end{enumerate}
In topology, one can not iterate this process very far before things become significantly more complicated. The only non-nilpotent elements to choose from are those in $\pi_0(\mathbb{S}^{\mathrm{top}})\cong \mathbb{Z}$, which are detected rationally. Therefore the Moore spectrum $\mathbb{S}^{\mathrm{top}}/x$ is of type 1, implying that any self-map in step (2) is forced to be a $v_1$-self map. In short, one is only able to find $v_1$-periodic elements in $\pi_*(\mathbb{S}^{\mathrm{top}})$ via this technique.

The situation becomes more delicate when one starts this program instead with a more general type $n$ spectrum $\mathrm{X}$. Such spectra exist and always have $v_n$-self-maps of some degree by the Periodicity Theorem \cite{DHS88,HS98}. However, constructing explicit models for type $n$ spectra is a very hard task, rendering the task of identifying $v_n$-periodic families in $\pi_*(\mathbb{S}^{\mathrm{top}})$ just as complicated. At height 2, there have been constructions of explicit type $2$ complexes. At other heights, there are currently no explicit models \cite{BehHilHopMah09_v2_32, BhaEggMah17_v2onA1, BhaEgg20_v2onZ}. 

As we have seen, there are many more non-nilpotent elements in $\pi_{*,*}^{\mathbb{C}}(\mathbb{S})$, and thus motivic spectra can hold multiple periodicities. Following Adams's program, Andrews showed the following.

\begin{thm}[\cite{Andrews18}]
    There is a nonnilpotent self-map
    \[
    w_1^4:\ss/\eta[20,12] \to \ss/\eta.
    \]
    Furthermore, let $\alpha \in \pi_{s,w}^{\mathbb{C}}(\mathbb{S})$ be one of the motivic stable elements $\{\nu, \nu^2, \nu^3, \overline{\sigma}, \overline{\sigma}\nu\}$,\footnote{It is insightful to compare these classes with those appearing in Adams's work by remembering their names on the Adams $\mathrm{E}_2$-page: $\eta$ is detected by $h_1$ and $\epsilon$ is detected by $c_0$; $\nu$ is detected by $h_2$ and $\overline{\sigma}$ is detected by $c_1$.} and let $w_1^{4k}(\alpha)$ be the composition
    \[
    w_1^{4k}(\alpha): \ss[s+2+20k, w+1+12k] \xrightarrow{\widetilde{\alpha}} \ss/\eta[20k, 12k] \xrightarrow{w_1^{4k}} \ss/\eta \xrightarrow{p} \ss[2,1].
    \]
    Then for any $k \geq 1$, $w_1^{4k}(\alpha) \in \pi^\mathbb{C}_{20k+s, 12k+w}(\mathbb{S})$ is essential.
\end{thm}

These are the first $w_1$-periodic classes to be identified in $\pi_{*,*}^\mathbb{C}(\mathbb{S})$. Analgous to how Adams's classes are 2-torsion $v_1$-periodic elements of $\pi_*(\mathbb{S}^{\mathrm{top}})$, Andrews's classes are $\eta$-torsion $w_1$-periodic elements of $\pi_{*,*}^\mathbb{C}(\mathbb{S})$. 

The more interesting part of the $v_1$-periodic portion of $\pi_*(\mathbb{S}^\mathrm{top})$ are the classes which are not simple 2-torsion. These classes exhibit interesting number theoretic properties in the particular degree of their 2-torsion.\footnote{By work of Friedlander \cite{Friedlander76}, for $p \equiv 3,5 \, (8)$ there is a map of spectra $\mathrm{GW}(\mathbb{F}_p) \to \mathbb{S}_{\mathrm{K}(1)}$ exhibiting the topological hermitian K-theory of $\mathbb{F}_p$ as the connective cover of the $\mathrm{K}(1)$-local sphere. In other words, there is an equivalence $\mathrm{GW}(\mathbb{F}_p)_{\mathrm{K}(1)} \simeq \mathbb{S}_{\mathrm{K}(1)}$.} In recent work of Isaksen, Kong, Li, Ruan, and Zhu, the analogous exotic story is studied \cite{IsaKonLiYuaZhu25}. The authors find new families in $\pi_{*,*}^\mathbb{C}(\mathbb{S})$ which are $w_1$-periodic and whose degree of $\eta$-torsion parallels the degree of $2$-torsion in the aforementioned $v_1$-periodic families which are not simple 2-torsion. The precise degree of torsion differs in a few notable places. For example, the topological classes $P^k(h_2^{\mathrm{top}})$ survives the $\textbf{ASS}(\mathbb{S}^{\mathrm{top}})$ to yield $v_1$-periodic classes $v_1^{4k}(\nu) \in \pi_{8k+3}(\mathbb{S}^{\mathrm{top}})$. The exotic analogue of these classes are $g^{4k}(h_3)$. Isaksen et al. show that these classes support differentials in the $\textbf{mASS}^\mathbb{C}(\mathbb{S})$ for all $k \geq 1$.\footnote{One can see the $\mathbb{R}$-motivic image-of-J in Belmont, Isaksen, and Kong's computation of $\pi_{*,*}^\mathbb{R}(\mathrm{L})$ \cite{belmontisaksenkong-v1R}. It is worth noting that the degree of 2-torsion higher than one might expect; for example, the third motivic Hopf map $\sigma$ is 32-torsion, whereas topologically it is 16-torsion. We thank J.D. Quigley for bringing this to our attention.} 


An interesting property of these $w_1$-periodic classes is that they survive under the map $\mathbb{S} \to \mathbb{S}/\tau$. By the identification of the $\mathbb{C}$-motivic homotopy of $\mathbb{S}/\tau$ with the $\mathrm{E}_2$-page of the topological Adams--Novikov spectral sequence, this allows one to see exotic periodicity manifest in the cohomology of $\euscr{M}_{\mathrm{fg}}$, the moduli stack of formal groups. Note that it has well been established that $\alpha_1$, which detects $\eta$, is non-nilpotent in this cohomology.

There has also been a construction of $w_1$-periodic elements in Quigley's work on the motivic Mahowald invariant for $F=\mathbb{C}$ and $\mathbb{R}$ \cite{Quigley19, Qui21_RC2mahowald}. This represents the only study of $w_1$-periodicity in $\mathrm{SH}(\mathbb{R})$.

As of the writing of this survey, there is currently no study of $w_1$-periodicity over fields other than $F=\mathbb{C}$ or $\mathbb{R}$ and no study of $w_n$-periodicity for $n \geq 2$ over any field.

\section{Future directions}
\label{section:future}
We conclude by highlighting some future directions and interesting problems.

\textbf{Problem:} \emph{Compute the Hurewicz image of $\mathrm{mmf}$ over $\mathbb{C}$.}

A useful way to detect $v_n$-periodicity in $\pi_*(\mathbb{S}^{\mathrm{top}})$ is by computing the Hurewicz image of chromatically interesting spectra. One example observed previously in this paper is the Hurewicz image of $\mathbb{S} \to \mathrm{ko}$, which consists of $\{1, v_1^{4k}(\eta), v_1^{4k}(\eta^2): k \geq 0\}$.

Many $v_2$-periodic families in $\pi_*(\mathbb{S}^{\mathrm{top}})$ were determined by work of Behrens, Mahowald, and Quigley, where they computed the Hurewicz image of $\mathbb{S} \to\mathrm{tmf}$ \cite{BehManQui23_tmfhurewiczimage}. We expect that all of these classes will lift to $\pi_{*,*}^\mathbb{C}(\mathbb{S})$ under the unit map $\mathbb{S} \to \mathrm{mmf}$; some of this has been confirmed by the work of Isaksen, Wang, and Xu \cite{IWX23}. Further study of this problem, in the context of $\mathrm{Syn}_{\mathrm{BP}}$, was conducted by Carrick and Davies \cite{cardav25_periodicfamilies}. Given that the $v_1$-periodicity of $\pi_{*,*}^\mathbb{C}(\mathbb{S})$ is very similar to the $v_1$-periodicity of $\pi_*(\mathbb{S}^{\mathrm{top}})$, with the only difference being the presence of $\tau$ and the non-nilpotence of $\eta$, we expect that the $v_2$-periodic elements of $\pi_{*,*}^\mathbb{C}(\mathbb{S})$ are not too far removed from the topological case.

\textbf{Problem:} \emph{Construct more ``motivic modular forms" spectra.}

There are many other topological spectra of modular forms (for example, see \cite{HilLaw16_tmflevelstructures}) which are of great utility. Many of these spectra are constructed by using the obstruction theory of Goerss, Hopkins, Miller, and Lurie, or alternatively via the techniques of spectral algebraic geometry \cite{Lurie09_ellipticsurvey, DouFraHenHil14_tmfbook}. A construction of motivic analogues of any of these spectra would provide very useful tools to study $\mathrm{SH}(F)$. As we have previously noted, for $F\neq \mathbb{C}$, there is as of yet no construction of even $\mathrm{mmf}$. Construction of such a spectrum would create great progress towards understanding how $v_2$-periodicity in $\pi_{*,*}^F(\mathbb{S})$ changes as the base field varies. 

\textbf{Problem:} \emph{Use the $\textup{\textbf{SliceSS}}^\mathbb{C}(\mathrm{mmf})$ to study $\Pi_k^\mathbb{C}(\mathbb{S})$.}

In their computation of $\Pi_1^F(\mathbb{S})$ and $\Pi_2^F(\mathbb{S})$, Rondigs, Spitzweck, and \O stv\ae r determine many differentials in the $\textbf{SliceSS}^F(\mathbb{S})$ by first studying the $\textbf{SliceSS}^F(\mathrm{kq})$ and then pulling back data along the unit map $\mathbb{S} \to \mathrm{kq}$ \cite{RSOfirst, RSOsecond}. As $\mathrm{mmf}$ is a better approximation to the sphere than $\mathrm{kq}$, we expect that lifting these techniques and pulling back slice data along the unit map $\mathbb{S} \to \mathrm{mmf}$ will lead to a better understanding of Milnor--Witt stems of higher degree. Moreover, the coefficient ring $\pi_{*,*}^\mathbb{C}(\mathrm{mmf})$ is completely understood (see \cite{isaksen_ASS_mmf, IKLRZ24_ANSSmmf}, for example), so that the actual determination of the differentials in the $\textbf{SliceSS}^\mathbb{C}(\mathrm{mmf})$ should not be too complicated. This would be especially useful if one were to construct a motivic modular forms spectrum which is stable under base change.

\textbf{Problem:} \emph{Construct resolutions of a motivic analogue of the $\mathrm{K}(2)$-local sphere.}

Another entry point to studying $v_n$-periodicity is by constructing a \emph{finite resolution} of the $\mathrm{K}(n)$-local sphere. For example, we saw previously a resolution of the $\mathrm{K}(1)$-local sphere
\[
\mathbb{S}^{\mathrm{top}}_{\mathrm{K}(1)} \to \mathrm{KO} \xrightarrow{\psi^3-1} \mathrm{KO},
\]
which has a motivic analogue by the recent work of Balderrama, Ormsby, and Quigley \cite{BOQ23}.

There are many resolutions of the $\mathrm{K}(2)$-local sphere at the prime 2 \cite{BobGoe18_Topologicalresolutions, Henn19_centralizer, BeaBobHen25}. These resolutions involve various fixed points of the Lubin--Tate spectrum $\mathrm{E}_2^{hG}$, where $G \subset \mathbb{G}_2$ is a finite subgroup of the height 2 Morava stabilizer group. One expects that, just as in the height 1 case, there are lifts of $\mathbb{S}^{\mathrm{top}}_{\mathrm{K}(2)}$ and its resolutions to $\mathrm{SH}(F)$. A naive candidate for a motivic analogue of the $\mathrm{K}(2)$-local sphere is the localization at the motivic spectrum $\mathrm{BPGL}\langle 2 \rangle/(2, v_1) = \mathrm{K}(2)^{\mathrm{mot}}$. By the motivic Landweber exact functor theorem of Naumann, Spitzweck, and \O stv\ae r \cite{motivic-landweber}, there is a motivic Lubin--Tate spectrum $\mathrm{E}_2^{\mathrm{mot}}$ over any base field. Additionally, Mazel--Gee has shown that the Morava stabilizer group $\mathbb{G}_n$ also acts on the motivic Lubin--Tate spectrum \cite{motivic-E-theory}.\footnote{However, it is not clear that $\mathbb{G}_n$ represents \emph{all} of the $\mathbb{E}_\infty$-automorphisms of $\mathrm{E}_n^{\mathrm{mot}}$.} Understanding resolutions of $\mathbb{S}_{\mathrm{K}(2)^{\mathrm{mot}}}$ would give a new perspective on $v_2$-periodicity in stable motivic homotopy theory.

\textbf{Problem:} \emph{Identify more $w_1$-periodic classes in $\pi_{*,*}^\mathbb{C}(\mathbb{S}).$}

The work of Andrews and of Isaksen, Kong, Li, Ruan, and Zhu has identified many classes in $\pi_{*,*}^\mathbb{C}(\mathbb{S})$ as $w_1$-periodic \cite{Andrews18, IsaKonLiYuaZhu25}. In Andrews's work, he constructs a $w_1^4$-self-map of $\mathbb{S}/\eta$ and imitates techniques of Adams techniques for producing $v_1$-periodic phenomena in $\pi_*(\mathbb{S}^{\mathrm{top}})$ via the $v_1^4$-self-map of $\mathbb{S}^{\mathrm{top}}/2$. In Isaksen et al.'s work, they identify classes in the $\mathrm{E}_2$-page of the $\textbf{mASS}^\mathbb{C}(\mathbb{S})$ which are not simple $h_1$-torsion and survive to give $w_1$-periodic elements in $\pi_{*,*}^\mathbb{C}(\mathbb{S})$ which are not simple $\eta$-torsion.

There is much more $w_1$-periodic phenomena which appears in the $\mathrm{E}_2$-page of the $\textbf{mASS}^\mathbb{C}(\mathbb{S})$. A systematic study of these classes, especially those which survive to $\pi_{*,*}^\mathbb{C}(\mathbb{S})$, could bring us closer to a conceptual understanding of $w_1$-periodicity. 

\textbf{Problem:} \emph{Compute the homotopy groups $\pi_{*,*}^\mathbb{C}(w_1^{-1}(\mathbb{S}/\eta)).$}

Consider the telescope 
\[
\mathrm{T}(1) = \mathrm{colim}(\mathbb{S}^{\mathrm{top}}/2 \xrightarrow{v_1^4}\mathbb{S}^{\mathrm{top}}/2[-8] \xrightarrow{v_1^4}\cdots).
\]
An important computation in the proof of the height 1 telescope conjecture is that of the homotopy groups $\pi_{*}(\mathrm{T}(1))$ \cite{Mah81}. In the same way, we can consider the telescope 
\[
\mathrm{T}(w_1) = \mathrm{colim}(\mathbb{S}/\eta \xrightarrow{w_1^4}\ss/\eta[-20,-12] \xrightarrow{w_1^4}\cdots).
\]
The computation of the homotopy groups $\pi_{*,*}^\mathbb{C}(\mathrm{T}(w_1))$ will be of vital interest in studying a $\mathbb{C}$-motivic reformulation of the telescope conjecture.

\textbf{Problem} \emph{Study $w_2$-periodic phenomena.}

The Periodicity Theorem of Hopkins and Smith implies that if $\mathrm{X} \in \mathrm{Sp}^\omega$, then there is a $v_n$-self-map $v:\mathrm{X}[d] \to \mathrm{X}$ inducing an isomorphism after smashing with $\mathrm{K}(n)$ \cite{HS98}. As a consequence, the cofiber $\mathrm{X}/v$ admits a $v_{n+1}$-self-map. These self-maps are useful for detecting periodicity in homotopy groups. In practice, it is difficult to identify explicit $v_n$-self-maps. A technique that can be useful is to start with the sphere $\mathbb{S}$ and iteratively take cofibers by self-maps, then use the Adams spectral sequence to determine the particular degree of a $v_{n+1}$-self-map. 

For example, at the prime 2, multiplication by 2 gives a $v_0$-self-map
\[
\mathbb{S}^{\mathrm{top}} \xrightarrow{2}\mathbb{S}^{\mathrm{top}} \to \mathbb{S}^{\mathrm{top}}/2.
\]
As we have seen, Adams showed that the cofiber admits a $v_1^4$-self-map
\[
\mathbb{S}^{\mathrm{top}}/2[8] \xrightarrow{v_1^4}\ss^{\mathrm{top}}/2 \to \ss^{\mathrm{top}}/(2, v_1^4)
\]
which allows on to identify many 8-periodic families of elements in $\pi_*(\mathbb{S}^{\mathrm{top}})$ \cite{Adams-JX-IV}. Work of Behrens, Hill, Hopkins, and Mahowald shows that the cofiber $\mathbb{S}/(2, v_1^4)$ admits a $v_2^{32}$-self-map
\[
\ss^{\mathrm{top}}/(2, v_1^4)[192] \xrightarrow{v_2^{32}} \ss^{\mathrm{top}}/(2, v_1^4) \to \ss^{\mathrm{top}}/(2, v_1^4, v_2^{32}),
\]
which allows one to identify many $192$-periodic families of elements in $\pi_*(\mathbb{S}^{\mathrm{top}})$ \cite{BehHilHopMah09_v2_32}.\footnote{It was initially claimed by Davis and Mahowald that $\ss^{\mathrm{top}}/(2, v_1^4)$ admitted a $v_2^8$-self-map \cite{DavMah81_v1v2periodicity}, which is not the case. However, it is conjectured that there are other finite complexes which admit $v_2^8$-self-maps, allowing one to identify 48-periodic families of elements in $\pi_*(\mathbb{S}^{\mathrm{top}})$.} It is currently unknown the degree of a minimal $v_3$-self-map on $\mathbb{S}^{\mathrm{top}}/(2, v_1^4, v_2^{32})$, although one must exist by the Periodicity Theorem.

Work of Krause guarantees the existence of $\mathbb{C}$-motivic exotic self-maps \cite{Kra18}, and it is expected that the above pattern has an exotic motivic analogue. To be precise, it is conjectured by Andrews that the cofiber of the $w_1^4$-self-map on $\ss/\eta$ admits a $w_2^{32}$-self-map
\[
\ss/(\eta, w_1^4) [416,224]\xrightarrow{w_2^{32}}\ss/(\eta, w_1^4) \to \ss/(\eta, w_1^4, w_2^{32}).
\]
This would identify many $(416, 224)$-periodic families of elements in $\pi_{*,*}^\mathbb{C}(\mathbb{S})$, and would be the first identification of $w_2$-periodic elements of any kind. As the computation of motivic stable stems is yet to enter this range, it would be much simpler to study $w_2$-periodicity if one could find a finite spectrum with a $w_2$-self-map of smaller period. 

Even still, the first element in such a $(416, 224)$-periodic familiy may potentially be in a computable range. We make the following observation. The class $\overline{\kappa} \in \pi_{20}(\mathbb{S}^{\mathrm{top}})$ is $v_2^{32}$-periodic. Under the Chow degree 0 isomorphism, the class $g \in \mathrm{Ext}^{20,4}_{\aeu^\vee_{\mathrm{top}}}(\mathbb{S}^{\mathrm{top}})$, which detects $\overline{\kappa}$, corresponds to the element $g_2 \in \mathrm{Ext}_{\aeu^\vee_\mathbb{C}}^{44,4,24}(\mathbb{S})$. This class survives the $\textbf{mASS}^\mathbb{C}(\mathbb{S})$ to give an element $\overline{\kappa}_2 \in \pi_{44,24}^\mathbb{C}(\mathbb{S})$.\footnote{Note that $\overline{\kappa}_2$ is non-nilpotent \cite[Example 6.18]{Kra18}} We conjecture that $\overline{\kappa}_2$ is $w_2^{32}$-periodic.\footnote{It is conjectured by Mahowald that $\nu \in \pi_3(\mathbb{S}^{\mathrm{top}})$ is $v_2^8$-periodic and $\epsilon \in \pi_8(\mathbb{S}^{\mathrm{top}})$ is $v_2^{16}$-periodic \cite[Chapter 13]{DouFraHenHil14_tmfbook}. If this is shown to be true, then we conjecture that $\sigma \in \pi_{7,4}^\mathbb{C}(\mathbb{S}^{\mathrm{top}})$ is $w_2^8$-periodic and $\overline{\sigma} \in \pi_{19, 11}^\mathbb{C}(\mathbb{S})$ is $w_2^{16}$-periodic.}

\textbf{Problem:} \emph{Study the interaction between $v_n$-periodicity and $w_n$-periodicity.}

Classes in $\pi_{*,*}^\mathbb{C}(\mathbb{S})$ may posses multiple periodicities. As $\eta$ is $v_1$-periodic, every $\eta$-periodic class is $v_1$-periodic. There are also classes in the Hurewicz image of $\mathrm{mmf}$ which are $v_2$-periodic, such as $\bar{\kappa}$, which is also $w_1^4$-periodic.

However, there is interesting interaction which is dependent on choice of base field. For example, we saw that $\pi_{*,*}^\mathbb{C}(\eta^{-1}\mathbb{S})$ is entirely 2-torsion, hence $(\eta^{-1}\mathbb{S})_{\mathbb{Q}} \simeq 0$ over $\mathbb{C}$. However, we saw that there is interesting torsion-free phenomena in $\pi_{*,*}^\mathbb{R}(\eta^{-1}\mathbb{S})$ stemming from the isomorphism of $\mathrm{W}(\mathbb{R}) \cong \mathbb{Z}$. This implies that $(\eta^{-1}\mathbb{S})_{\mathbb{Q}} \not\simeq 0$ over $\mathbb{R}$, and a similar result holds over all fields $F$ such that $\mathrm{W}(F)$ is not entirely torsion.

A particularly interesting question is: if one were to construct $w_1$-periodic families in $\pi_{*,*}^F(\mathbb{S})$, how do these families behave upon $\mathrm{KGL}/2$-localization, and how does the answer vary as the base field changes?

\textbf{Problem:} \emph{Study odd primary periodicity.}
Most of the computational work in motivic homotopy theory has been at the prime 2. However, there is interesting phenomena that occurs at odd primes \cite{Sta23, MorPetTat-BPGL1}. While $\eta$ is 2-torsion, hence does not appear at odd primes, there are other elements of $\pi_{*,*}^F(\mathbb{S})$ which are visible to odd primes which lift topological elements and are non-nilpotent.

For example, at the prime 3, the element $\beta_1 \in \pi_{10,6}^\mathbb{C}(\mathbb{S})$ is non-nilpotent. Therefore, one may ask about $\beta_1$-periodic phenomena.\footnote{At the prime 3, the $\textbf{mANSS}^\mathbb{C}(\beta_1^{-1}\mathbb{S})$ has $\mathrm{E}_2$-page given by $\beta_1^{-1}\mathrm{Ext}_{\mathrm{BP}_*\mathrm{BP}}(\mathbb{S}^{\mathrm{top}}) \otimes \mathbb{Z}[\tau]$. The $\beta_1$-localized Ext group appearing here has not been fully computed, but was thoroughly studied by Belmont \cite{Belmont20_localizingExt}.} One would assume that, just as there is a non-nilpotent self-map of $\mathbb{S}/\eta$ giving access to higher exotic periodicity, one should also be able to find a non-nilpotent self-map of $\mathbb{S}/\beta_1$. Just as Isaksen et al. push exotic periodicity to the cofiber of $\tau$ to discover patterns in the cohomology of $\euscr{M}_{\mathrm{fg}}$, one would also assume that odd primary exotic periodicity will also manifest in the odd primary topological Adams--Novikov $\mathrm{E}_2$-page.

Notice that $\eta$ is $v_1$-periodic at the prime 2, $w_1^4$ is $v_2$-periodic at the prime 2, and $\beta_1$ is $v_2$-periodic at the prime 3. Work in progress of Burklund, Levy, and Pstragowski suggests that, at least over $\mathbb{C}$, at the prime $p$, one should expect to see exotic periodicity stemming from a $v_h$-periodic element, where $h = n(p-1)$. For example, at the prime 5, one should expect to find a $v_4$-periodic element in $\pi_{*,*}^\mathbb{C}(\mathbb{S})$ which is non-nilpotent, allowing one to run this exotic periodicity game.\footnote{For primes $p>3$, there is some handle on height $p-1$-phenomena due to Ravenel's solution to the $p>3$-primary Kervaire invariant problem \cite{Rav78_oddprimarykervaire}.}

\textbf{Problem:} \emph{Determine the structure of the Balmer spectrum $\mathrm{Spc}(\mathrm{SH}(F)^{\mathrm{cell}})$.}

For an idempotent-complete stably symmetric monoidal $\infty$-category $\euscr{C}$, we denote by $\mathrm{Spc}(\euscr{C})$ the \emph{Balmer spectrum} of $\euscr{C}$ \cite{AokBarCheSchSte25_higherzariski2ring}. This is a topological space which is defined by the support data of $\euscr{C}$ in the same way the Zariski spectrum $\mathrm{Spec}(R)$ of a commutative ring $R \in \mathrm{CRing(\mathrm{Ab})}$ is defined by the prime ideals of $R$.

The Thick Subcategory Theorem of Hopkins and Smith may be reinterpreted as a computation of the Balmer spectrum of compact spectra $\mathrm{Sp}^\omega$ \cite{HS98}. The points of $\mathrm{Spc}(\mathrm{Sp}^\omega)$ correspond to the subcategories which are $\mathrm{K}(n)$-acyclic, meaning $\mathrm{K}(n) \otimes \mathrm{X} \simeq 0$, and the topology is determined by the fact that $\mathrm{K}(n) \otimes \mathrm{X}\simeq 0$ implies that $\mathrm{K}(n-1)\otimes \mathrm{X} \simeq 0$.

There has been some study of the motivic Balmer spectrum $\mathrm{Spc}(\mathrm{SH}(F)^\omega)$. Work of Heller and Ormsby shows that there is a surjective map
\[
\rho_\bullet:\mathrm{Spc}(\mathrm{SH}(F)^\omega) \to \mathrm{Spec}^h(\mathrm{K}^{\mathrm{MW}}_*(F))
\]
to the homogeneous spectrum of Milnor--Witt K-theory \cite{HelOrm18_ttgSHF}. By Morel's isomorphism $\Pi_0^F(\mathbb{S}) \cong \mathrm{K}^{\mathrm{MW}}_*(\mathrm{F})$ one can think of this as the motivic analogue of the map
\[
\mathrm{Spc}(\mathrm{Sp}^\omega) \to \mathrm{Spec}(\pi_0(\mathbb{S}^{\mathrm{top}})) \cong \mathrm{Spec}(\mathbb{Z}).
\]
The homogeneous spectrum $\mathrm{Spec}^h(\mathrm{K}^{\mathrm{MW}}_*(F))$ has been computed for all fields of characteristic not 2 by Thornton \cite{Thornton16_speckmw}. Thus, describing the Balmer spectrum of motivic spectra amounts to understanding the fibers of $\rho_\bullet$.

This is a quite complicated task. The topology of the Balmer spectrum is related to the ways in which periodic phenomena interact with each other. As we have seen, this is not well understood. Finite motivic spectra may posess multiple periodicities which glue together in sophisticated ways. For example, we have seen that the sphere spectrum $\mathbb{S}$ admits $p$-periodicity for any prime $p$, $\eta$-periodicity, and $\beta_1$-periodicity. The relationship between these periodicities will inform the structure of $\mathrm{Spc}(\mathrm{SH}(F)^{\omega})$. There has been some work in understanding this structure when $F=\mathbb{C}$ due to Krause, Hornbostel, and Joachimi \cite{Kra18,Horn18_motivicnilpotence, Joachimi20}. 

In related work, Du and Vishik have studied the Balmer spectrum of the isotropic stable motivic homotopy category \cite{DuVis25_balmerspectrumSHKiso}. Their work has also been influential to that of Balmer and Gallauer \cite{BalGal25_specArtinMotives, Gal25_periodicsmotivic}. Finally, there has been substantial progress in understanding the Balmer spectrum of Voevodsksy's derived category of motives (for example, see \cite{Pet13_primesinSpcArtintateMotives, Gal19_ttgTateMotives}).

\textbf{Problem:} \emph{Find a way to effectively detect nilpotence.}

The Nilpotence Theorem of Devinatz, Hopkins, and Smith is an important structural result in chromatic homotopy theory \cite{DHS88}. One particular form is that if $\mathrm{R} \in \mathrm{Sp}$ is a homotopy ring spectrum, then $x \in \pi_*(\mathrm{R})$ is nilpotent if and only if the $\mathrm{MU}$-Hurewicz image $\mathrm{MU}_*(x) \in \mathrm{MU}_*(\mathrm{R})$ is nilpotent. In other words, $\mathrm{MU}$ detects nilpotence. A related result is the May Nilpotence conjecture, now proven by Mathew, Naumann, and Noel, which states that if $\mathrm{R} \in \mathrm{CAlg}(\mathrm{Sp})$, then $\mathrm{MU} \otimes \mathrm{R}\simeq 0$ if and only if $\mathrm{H}\mathbb{Z}\otimes \mathrm{R} \simeq 0$ \cite{MatNauNoe15_MayNilpotenceConjecture}. 

From many perspectives, the algebraic cobordism spectrum $\mathrm{MGL}$ serves as a motivic analogue of the complex cobordism spectrum $\mathrm{MU}$. For example, if $\mathrm{E}$ is an oriented motivic cohomology theory, then it receives a ring map $\mathrm{MGL} \to \mathrm{E}$, and pushing forward the universal formal group law over $\mathrm{MGL}$ defines a formal group law over $\mathrm{E}$. 

Bachmann and Hahn have proved a version the motivic May Nilpotence conjecture, which states that $\mathrm{MGL} \otimes \mathrm{R}\simeq 0$ if and only if $\mathrm{H}\mathbb{Z} \otimes \mathrm{R}\simeq 0$ when $\mathrm{R}$ is a \emph{normed}\footnote{A normed motivic ring spectrum is an enhancement of being $\mathbb{E}_\infty$ and is a sufficient condition to satsify the motivic May nilpotence conjecture. See the work of Bachmann and Hoyois for details on normed motivic ring spectra \cite{BachmannHoyois-Norms}.} motivic ring spectrum \cite{BacHah22_MGLnilpotence}. However, $\mathrm{MGL}$ does not in general detect nilpotence. A simple example is the Hopf map $\eta \in \pi_{1,1}^F(\mathbb{S})$ which is non-nilpotent but not detected by $\mathrm{MGL}$. 

There are other spectra which detect nilpotent elements. For example, the motivic symplectic cobordism spectrum $\mathrm{MSp}$ detects $\eta$ \cite{PanWal23_MSLMSp}. Botvinik and Kochman have shown that the homotopy of the topological $\mathrm{MSp}$ is very complicated, containing $2^k$-torsion elements for all $k$ \cite{Koc93_MspandStableStems,BotKoc94_MSptorsion}. 
Consequentially, studying the motivic $\mathrm{MSp}$ may not seem very amenable to computation.
However, a surprising result of Bachmann and Hopkins shows that after inverting $\eta$, the algebraic symplectic cobordism spectrum is significantly less complicated: there is an isomorphism \cite{BH21}
\[
\pi^F_{*,*}(\eta^{-1}\mathrm{MSp}) \cong \mathrm{W}(F)[\eta^{\pm 1}, y_1, y_2, \dots],
\]
where $|y_i| = (2i,0)$. We imagine that other localizations of $\mathrm{MSp}$ may be useful for detecting exotic nilpotence and periodicity. See also the work of Hornbostel and of Bachmann and Hahn for robust discussions of nilpotence in $\mathrm{SH}(F)$ \cite{Horn18_motivicnilpotence,BacHah22_MGLnilpotence}.

\textbf{Problem:} \emph{Understand the impact of symplectic orientations on exotic chromatic structures.}

An interesting feature of $\mathrm{MSp}$ is that it is the universal motivic cohomology theory which is symplectically oriented, in that any symplectically oriented cohomology theory $\mathrm{E}$ receives a map $\mathrm{MSp} \to \mathrm{E}$.\footnote{Pulling back along the ring map $\mathrm{MSp} \to \mathrm{MGL}$ shows that every $\mathrm{GL}$-oriented motivic spectrum also admits a $\mathrm{Sp}$-orientation.} Let $\euscr{W}$ denote the so-called \emph{Walter ring}, which is the ring over which the universal formal ternary law is defined, where a formal ternary law is a three-term analogue of a formal group law \cite{DF23}. The observation that the tensor product of any 3 symplectic bundles is again symplectic leads to a formal ternary law associated to any symplectically oriented motivic spectrum. This implies that there is a map $\euscr{W} \to \mathrm{E}_{2*,*}$ which sends the universal formal ternary law over $\euscr{W}$ to one over $\mathrm{E}_{2*,*}$.

Just as there is an isomorphism $\euscr{L} \xrightarrow{\simeq} \mathrm{MGL}_{2*, *}$ between the Lazard ring and a subring of the universal oriented cohomology theory, it is conjectured by Coulette, D\'eglisse, Fasel, and Hornbostel that the classifying map $\euscr{W} \to \mathrm{MSp}_{2*, *}$ is an isomorphism \cite{ColDegFasHor24_2groupsandFTLs}. As is paralleled by the complexity of the homotopy of $\mathrm{MSp}$, a complete description of $\euscr{W}$ is not available.
Very recent work of Huang has investigated the role of formal ternary laws and $\mathrm{MSp}$ in understanding symplectic orientations in the $\eta$-periodic category $\eta^{-1}\mathrm{SH}(F)$ \cite{Huang26_etaperiodicFTL}. A particularly interesting result is the construction of an $\eta$-periodic Walter ring $\euscr{W}^\eta$ and isomorphisms 
\[
\euscr{L}[1/2] \xrightarrow{\simeq} \euscr{W}^{\eta}[1/2] \xrightarrow{\simeq} \pi_{2*, *}^F(\eta^{-1}\mathrm{MSp}[1/2]),
\]
at least when $\mathrm{W}(F) \cong \mathbb{Z}$.

The utility of complex-orientable cohomology theories in studying chromatic homotopy theory inclines one to believe that the study of symplectic orientations could play an important role in understanding chromatic motivic homotopy theory. A particularly interesting question is the following. A useful way to construct complex-oriented spectra is by way of the Landweber exact functor theorem. A motivic analogue of this theorem, which produces $\mathrm{GL}$-oriented motivic spectra, has been proven by Naumann, Spitzweck, and \O stv\ae r \cite{motivic-landweber}. This leads one to believe that there should be a symplectic Landweber exact functor theorem which allows one to construct $\mathrm{Sp}$-oriented motivic spectra.











\printbibliography

\end{document}